\documentclass{amsart}
\usepackage{amssymb}
\usepackage{mathrsfs}
\usepackage{stmaryrd} 
\usepackage{bbm}
\usepackage{oldgerm}
\usepackage[francais]{babel}
\usepackage[T1]{fontenc}
\usepackage[latin1]{inputenc}
\usepackage[all]{xy}
\usepackage{hyperref}

\newtheorem{teo}[subsection]{Théorème}
\newtheorem{prop}[subsection]{Proposition}
\newtheorem{cor}[subsection]{Corollaire}

\theoremstyle{definition}

\newtheorem{defi}[subsection]{Définition}
\newtheorem{rema}[subsection]{Remarque}

\numberwithin{equation}{subsection}

\newcommand{\mH}{{\mathbb H}}

\newcommand{\mQ}{{\mathbb Q}}
\newcommand{\mN}{{\mathbb N}}

\newcommand{\mZ}{{\mathbb Z}}

\newcommand{\mG}{{\mathbb G}}
\newcommand{\mK}{{\mathbb K}}

\newcommand{\mV}{{\mathbb V}}

\newcommand{\bB}{{\bf B}}
\newcommand{\bD}{{\bf D}}
\newcommand{\bL}{{\bf L}}
\newcommand{\bT}{{\bf T}}

\newcommand{\Et}{{\bf \acute{E}t}}

\newcommand{\Ens}{{\bf Ens}}

\newcommand{\bMod}{{\bf Mod}}
\newcommand{\MOD}{{\rm MOD}}
\newcommand{\bIH}{{\bf IH}}

\newcommand{\et}{{\rm \acute{e}t}}
\newcommand{\fet}{{\rm f\acute{e}t}}

\newcommand{\ad}{{\rm ad}}

\newcommand{\zar}{{\rm zar}}

\newcommand{\coh}{{\rm coh}}

\newcommand{\Dolb}{{\rm Dolb}}
\newcommand{\sol}{{\rm sol}}

\newcommand{\atf}{{\rm atf}}

\newcommand{\rf}{{\rm f}}

\newcommand{\HT}{{\rm HT}}

\newcommand{\Spec}{{\rm Spec}}

\newcommand{\Spf}{{\rm Spf}}

\newcommand{\ob}{{\rm Ob}}

\newcommand{\tot}{{\rm tot}}

\newcommand{\cont}{{\rm cont}}

\newcommand{\id}{{\rm id}}

\newcommand{\Sym}{{\rm Sym}}

\newcommand{\Hom}{{\rm Hom}}

\newcommand{\End}{{\rm End}}

\newcommand{\bMH}{{\bf MH}}

\newcommand{\rC}{{\rm C}}
\newcommand{\rD}{{\rm D}}

\newcommand{\rH}{{\rm H}}
\newcommand{\rT}{{\rm T}}

\newcommand{\rR}{{\rm R}}
\newcommand{\rS}{{\rm S}}

\newcommand{\rW}{{\rm W}}

\newcommand{\oK}{{\overline{K}}}

\newcommand{\oR}{{\overline{R}}}
\newcommand{\oS}{{\overline{S}}}

\newcommand{\oU}{{\overline{U}}}

\newcommand{\oX}{{\overline{X}}}
\newcommand{\oY}{{\overline{Y}}}

\newcommand{\ogg}{{\overline{g}}}

\newcommand{\ox}{{\overline{x}}}
\newcommand{\oy}{{\overline{y}}}

\newcommand{\oeta}{{\overline{\eta}}}

\newcommand{\ocB}{{\overline{\cB}}}

\newcommand{\upp}{{\underline{p}}}

\newcommand{\hA}{{\widehat{A}}}

\newcommand{\hRun}{{\widehat{R_1}}}

\newcommand{\halpha}{\widehat{\alpha}}

\newcommand{\cA}{{\mathscr A}}
\newcommand{\cB}{{\mathscr B}}
\newcommand{\cC}{{\mathscr C}}

\newcommand{\cE}{{\mathscr E}}
\newcommand{\cF}{{\mathscr F}}
\newcommand{\cG}{{\mathscr G}}

\newcommand{\cK}{{\mathscr K}}
\newcommand{\cL}{{\mathscr L}}

\newcommand{\co}{{\mathscr O}}
\newcommand{\cR}{{\mathscr R}}
\newcommand{\cS}{{\mathscr S}}

\newcommand{\cH}{{\mathscr H}}
\newcommand{\cM}{{\mathscr M}}
\newcommand{\cN}{{\mathscr N}}

\newcommand{\cU}{{\mathscr U}}
\newcommand{\cV}{{\mathscr V}}

\newcommand{\fN}{{\mathfrak N}}

\newcommand{\fS}{{\mathfrak S}}
\newcommand{\fV}{{\mathfrak V}}
\newcommand{\fX}{{\mathfrak X}}

\newcommand{\fgg}{{\mathfrak g}}

\newcommand{\tta}{{\tt a}}

\newcommand{\hM}{{\widehat{M}}}

\newcommand{\hR}{{\widehat{R}}}

\newcommand{\hY}{{\widehat{Y}}}
\newcommand{\hoR}{{\widehat{\oR}}}

\newcommand{\hcC}{{\widehat{\cC}}}

\newcommand{\bvd}{{\breve{d}}}
\newcommand{\bvu}{{\breve{u}}}

\newcommand{\bvoS}{{\breve{\oS}}}

\newcommand{\bvoX}{{\breve{\oX}}}
\newcommand{\bvocB}{{\breve{\ocB}}}

\newcommand{\bvcC}{{\breve{\cC}}}

\newcommand{\bvcF}{{\breve{\cF}}}

\newcommand{\bvnu}{{\breve{\nu}}}
\newcommand{\bvalpha}{{\breve{\alpha}}}
\newcommand{\bvsigma}{{\breve{\sigma}}}

\newcommand{\bvPhi}{{\breve{\Phi}}}

\newcommand{\coS}{{\check{\oS}}}

\newcommand{\coX}{{\check{\oX}}}
\newcommand{\coY}{{\check{\oY}}}
\newcommand{\cog}{{\check{\ogg}}}

\newcommand{\tE}{{\widetilde{E}}}

\newcommand{\tU}{{\widetilde{U}}}

\newcommand{\tX}{{\widetilde{X}}}
\newcommand{\tY}{{\widetilde{Y}}}

\newcommand{\tg}{{\widetilde{g}}}

\begin{document}

\title[Sur la correspondance de Simpson $p$-adique. 0]
{Sur la correspondance de Simpson $p$-adique.\\ 0~: une vue d'ensemble}
\author{Ahmed Abbes et Michel Gros}
\address{A.A. Laboratoire Alexander Grothendieck, ERL 9216 du CNRS, Institut des Hautes \'Etudes Scientifiques, 35 route de Chartres, 91440 Bures-sur-Yvette, France}
\address{M.G. CNRS UMR 6625, IRMAR, Université de Rennes 1,
Campus de Beaulieu, 35042 Rennes cedex, France}
\email{abbes@ihes.fr}
\email{michel.gros@univ-rennes1.fr}

\maketitle

\setcounter{tocdepth}{1}
\tableofcontents

\section{Introduction}

\subsection{}\label{higgs0-intro1}
En 1965, généralisant et précisant un résultat de Weil, Narasimhan et Seshadri \cite{ns} 
établissaient une correspondance bijective entre l'ensemble des classes d'équivalence de représentations 
unitaires irréductibles du groupe fondamental d'une surface de Riemann compacte $X$ de genre $\geq 2$, 
et l'ensemble des classes d'isomorphisme de fibrés vectoriels stables de degré $0$ sur $X$. 
La correspondance fut ensuite étendue pour toute variété complexe projective et lisse par Donaldson \cite{donaldson}. 
L'analogue pour les représentations linéaires quelconques est dû à Simpson~; 
pour obtenir une correspondance du même type que celle de Narasimhan et Seshadri, on a besoin 
d'ajouter au fibré vectoriel une structure supplémentaire. C'est la notion de fibré de Higgs 
qui fut d'abord introduite par Hitchin pour les courbes algébriques. 
Si $X$ est un schéma lisse sur un corps $K$, 
un {\em fibré de Higgs}  sur $X$ est un couple $(M,\theta)$ formé 
d'un $\co_X$-module localement libre de type fini $M$ et d'un morphisme $\co_X$-linéaire 
$\theta\colon M\rightarrow M\otimes_{\co_X}\Omega^1_{X/K}$
tel que $\theta\wedge \theta=0$. 
Le résultat principal de Simpson \cite{simpson1,simpson4,simpson2,simpson3} établit une équivalence de catégories
entre la catégorie des représentations linéaires de dimension finie (à valeurs complexes) 
du groupe fondamental d'une variété complexe projective et lisse
et celle des fibrés de Higgs semi-stables de classes de Chern nulles (cf. \cite{lepotier}).

\subsection{}\label{higgs0-intro2}
Les résultats de Simpson et les développements qu'ils ont suscités ont amené 
depuis quelques années la recherche d'un  analogue $p$-adique. 
Les premiers exemples d'une telle construction (qui n'utilisait pas encore la terminologie de fibrés de Higgs)
se trouvent rétrospectivement chez Hyodo \cite{hyodo1} qui avait traité le cas conceptuellement important des variations 
de structures de Hodge $p$-adiques, baptisées {\em systèmes locaux de Hodge-Tate}. 
L'approche la plus avancée au stade actuel est due à Faltings \cite{faltings3}. 
Elle généralise des résultats antérieurs de Tate, Sen et Fontaine, 
et repose sur sa théorie des extensions presque-étales \cite{faltings2}.
Elle s'inscrit dans le prolongement de ses travaux en théorie de Hodge $p$-adique, 
en particulier, ceux qui établissent l'existence de décompositions de Hodge-Tate \cite{faltings1}. Achevée, 
la {\em correspondance de Simpson $p$-adique} devrait ainsi naturellement fournir les meilleurs énoncés du type 
Hodge-Tate en théorie de Hodge $p$-adique. 
Mais au stade actuel, la construction de Faltings n'apparaît satisfaisante que pour les courbes, 
et même dans ce cas, de nombreuses questions, fondamentales, demeurent ouvertes.

\subsection{}\label{higgs0-intro3}
Nous développons dans ce travail \cite{ag1,ag3,ag2} une nouvelle 
approche de cette correspondance de Simpson $p$-adique, intimement liée à celle de Faltings, 
et inspirée du travail d'Ogus et Vologodsky \cite{ov} sur un analogue 
en caractéristique $p$ de la correspondance de Simpson complexe.
Le besoin de reprendre et de développer la construction de Faltings s'est fait sentir au regard
du nombre de résultats esquissés dans un article assez court et extrêmement dense \cite{faltings3}.
Signalons aussi que Tsuji a développé une autre approche pour la correspondance de 
Simpson $p$-adique \cite{tsuji3}. Par ailleurs, Deninger et Werner \cite{deninger-werner1} 
ont développé un analogue partiel de la théorie de Narasimhan et Seshadri pour les courbes $p$-adiques,
qui correspond aux fibrés de Higgs à champ nul dans la  correspondance de Simpson $p$-adique.  

\subsection{}\label{higgs0-intro4}
Soient $K$ un corps de valuation discrète complet de caractéristique
$0$, à corps résiduel parfait de caractéristique $p>0$, $\co_K$ l'anneau de valuation de $K$, 
$\oK$ une clôture algébrique de $K$, $\co_\oK$ la clôture intégrale de $\co_K$ dans $\oK$.
Soient $X$ un $\co_K$-schéma lisse de type fini à fibre géométrique générique $X_\oK$ intègre, 
$\ox$ un point géométrique de $X_\oK$, $\fX$ le schéma formel complété $p$-adique de $X\otimes_{\co_K}\co_{\oK}$.
Nous considérerons dans ce travail une situation logarithmique lisse plus générale (cf. \cite{ag1} 6.2 et \cite{ag3} 4.7).
Toutefois, pour simplifier la présentation, nous nous limitons dans cette introduction au cas lisse au sens usuel.   
On cherche à construire un foncteur de la catégorie des représentations $p$-adiques 
du groupe fondamental géométrique $\pi_1(X_\oK,\ox)$ ({\em i.e.}, les $\mQ_p$-représentations linéaires 
continues de dimension finie de $\pi_1(X_\oK,\ox)$)
dans la catégorie des $\co_\fX[\frac 1 p]$-fibrés de Higgs ({\em i.e.}, les couples $(\cM,\theta)$ formés 
d'un $\co_\fX[\frac 1 p]$-module localement projectif de type fini $\cM$ et d'un morphisme $\co_\fX[\frac 1 p]$-linéaire 
$\theta\colon \cM\rightarrow \cM\otimes_{\co_X}\Omega^1_{X/\co_K}$ tel que $\theta\wedge \theta=0$ \eqref{higgs0-md17}).
Suivant la stratégie de Faltings, qui n'est que partiellement achevée au stade actuel, 
ce foncteur devrait s'étendre à une catégorie strictement plus large que celle des représentations $p$-adiques de 
$\pi_1(X_\oK,\ox)$, appelée catégorie des {\em représentations généralisées} de $\pi_1(X_\oK,\ox)$.
Il serait alors une {\em équivalence de catégories} entre cette nouvelle catégorie et la catégorie des $\co_\fX[\frac 1 p]$-fibrés 
de Higgs. La motivation principale du présent travail est la construction d'une telle équivalence de catégories.
Lorsque $X_K$ est une courbe propre et lisse sur $K$, Faltings montre que les fibrés de Higgs associés aux 
``vraies'' représentations $p$-adiques de $\pi_1(X_\oK,\ox)$ 
sont semi-stables de pente nulle, et formule l'espoir que tous les fibrés de Higgs semi-stables de pente nulle
s'obtiennent ainsi. Cet énoncé correspondrait à la partie difficile du résultat de Simpson dans le cas complexe.

\subsection{}\label{higgs0-intro5}
La notion de représentation généralisée est due à Faltings.
Ce sont, en termes simplifiés, des représentations {\em semi-linéaires} $p$-adiques continues de $\pi_1(X_\oK,\ox)$ 
dans des modules sur un certain anneau $p$-adique muni d'une action continue de $\pi_1(X_\oK,\ox)$. 
L'approche de Faltings dans \cite{faltings3} pour construire 
un foncteur $\cH$ de la catégorie de ces représentations généralisées vers la catégorie des fibrés de Higgs,
comporte deux étapes.  Il définit d'abord $\cH$ pour les représentations généralisées 
qui sont $p$-adiquement proches de la représentation triviale, qu'il qualifie de {\em petites}.
Cette première étape est achevée en toute dimension. Dans la seconde étape, réalisée seulement 
pour les courbes, il étend le foncteur $\cH$ à toutes les représentations généralisées de 
$\pi_1(X_\oK,\ox)$ par descente. En effet, toute représentation
généralisée devient petite au-dessus d'un revêtement fini étale de $X_\oK$. 

Notre nouvelle approche, valable en toute dimension, permet de définir le foncteur $\cH$
sur la catégorie des représentations généralisées de $\pi_1(X_\oK,\ox)$ vérifiant une condition d'admissibilité à la Fontaine,
baptisées représentations généralisées {\em de Dolbeault}. Pour ce faire, nous introduisons des anneaux de périodes,
que nous appelons {\em algèbres de Higgs-Tate}, qui sont la principale nouveauté de notre approche par rapport à celle de Faltings.
Nous montrons que la condition d'admissibilité pour les coefficients rationnels correspond à la condition de petitesse de Faltings~;
mais elle est strictement plus générale pour les coefficients entiers. On notera que la construction de Faltings pour les  
coefficients rationnels petits est limitée aux courbes et qu'elle présente quelques difficultés que notre approche permet d'éviter.

\subsection{}\label{higgs0-intro6}
Nous procédons en deux étapes. Nous étudions d'abord dans \cite{ag1} le cas d'un schéma affine 
d'un type particulier, aussi qualifié par Faltings de {\em petit}. 
Nous abordons ensuite  dans \cite{ag3} les aspects globaux de la théorie.
La construction générale s'obtient à partir du cas affine par une technique de recollement recelant  
des difficultés inattendues. Nous utilisons pour ce faire le {\em topos de Faltings}, une variante fibrée de la notion de 
{\em topos co-évanescent} de Deligne que nous développons dans \cite{ag2}.

\subsection{}\label{higgs0-intro7}
Cette introduction présente, dans une situation géométrique simplifiée pour la clarté de l'exposition, un
résumé détaillé des grandes étapes conduisant à nos principaux résultats. Passons brièvement en revue
son contenu. On commence dans § \ref{higgs0-prg} par une courte digression sur les petites représentations généralisées 
dans le cas affine qui serviront comme intermédiaire pour l'étude des représentations de Dolbeault.  
La section \ref{higgs0-dth} résume la première partie \cite{ag1}. On introduit la notion de représentation généralisée de Dolbeault
pour un petit schéma affine et la notion compagnon de module de Higgs soluble, 
puis on construit une équivalence naturelle entre ces deux catégories.
On développe en fait deux variantes, une entière et une rationnelle plus subtile. 
On établit des liens entre ces notions et les conditions de petitesse de Faltings. 
On fait aussi le lien avec la théorie de Hyodo \cite{hyodo1}. 
La seconde partie \cite{ag3} est résumée dans les sections \ref{higgs0-rmd}
et \ref{higgs0-md}. Après une brève introduction au topos annelé de Faltings dans § \ref{higgs0-rmd},
on introduit les {\em algèbres de Higgs-Tate} \eqref{higgs0-rmd6}. 
La notion de {\em module de Dolbeault} qui globalise celle de représentation généralisée de Dolbeault, 
et la notion compagnon de {\em fibré de Higgs soluble} sont définies dans \ref{higgs0-md13}. 
Notre principal résultat \eqref{higgs0-md18} est l'équivalence de ces deux catégories. 
La preuve de ce résultat nécessite des énoncés d'acyclicité des algèbres de Higgs-Tate qu'on trouvera 
dans \ref{higgs0-md5} et \ref{higgs0-md7}, lesquels   permettent également de montrer la compatibilité de 
cette équivalence au passage à la cohomologie \eqref{higgs0-md181}. 
Nous étudions également la fonctorialité des diverses propriétés introduites par morphismes étales \eqref{higgs0-md20}, 
ainsi que leur caractère local pour la topologie étale (\ref{higgs0-md21}, \ref{higgs0-md22}, \ref{higgs0-md27}). 
Enfin, nous revenons dans cette situation globale  sur les liens logiques (\ref{higgs0-md24}, \ref{higgs0-md25}, \ref{higgs0-md26}),  
pour un fibré de Higgs, entre petitesse \eqref{higgs0-md23} et  caractère soluble.

Le lecteur trouvera au début des articles \cite{ag1,ag3} une description détaillée de leurs structures. 
L'article \cite{ag2}, étant d'un intérêt indépendant, a sa propre introduction.

\subsection*{Remerciements}
Ce travail n'aurait évidemment pu voir le jour sans les travaux de G. Faltings, 
en premier lieu celui sur la correspondance de Simpson $p$-adique \cite{faltings3}.
Nous lui exprimons notre profonde reconnaissance. 
La genèse de ce travail est immédiatement postérieure à un groupe de travail tenu à 
Rennes en 2008-09 sur son article \cite{faltings3}. Nous avons pu bénéficier à cette 
occasion du texte de l'exposé d'O. Brinon \cite{brinon1} et du travail de T. Tsuji \cite{tsuji5} 
présentant sa propre approche de la correspondance de Simpson $p$-adique. 
Ces deux textes nous ont été extrêmement utiles et nous sommes reconnaissants à leurs auteurs 
de les avoir spontanément mis à notre disposition. Nous remercions aussi O. Brinon, G. Faltings et T. Tsuji pour 
tous les échanges que nous avons eus avec eux sur des questions relatives à ce travail, 
et  A. Ogus pour les discussions éclairantes que nous avons eues avec lui sur son travail avec V. Vologodsky \cite{ov}. 
Le premier auteur (A.A.) remercie le Centre Émile Borel, l'Institut des Hautes Études Scientifiques
et l'Université de Tokyo pour leur hospitalité. Il remercie également 
les auditeurs d'un cours qu'il a donné sur ce sujet à l'Université de Tokyo durant l'automne
2010 et l'hiver 2011 dont les questions et remarques ont été précieuses pour mettre au point ce travail.
Le second auteur (M.G.) remercie l'Institut des Hautes Études Scientifiques
et l'Université de Tokyo pour leur hospitalité. 
Enfin, nous remercions  les participants de l'école d'été 
{\em Higgs bundles on $p$-adic curves and Representation theory} qui s'est tenue à Mainz en septembre 2012
et au cours de laquelle nos principaux résultats ont été exposés pour leurs remarques et leur intérêt stimulant.
Ce travail a bénéficié du soutien du programme ANR Théorie de Hodge $p$-adique 
et Développements (ThéHopaD) ANR-11-BS01-005.

\section{Notations et conventions}\label{higgs0-not}

{\em Tous les anneaux considérés dans cet article possèdent un élément unité~;
les homomorphismes d'anneaux sont toujours supposés transformer l'élément unité en l'élément unité.
Nous considérons surtout des anneaux commutatifs, et lorsque nous parlons d'anneau 
sans préciser, il est sous-entendu qu'il s'agit d'un anneau commutatif~; en particulier, 
il est sous-entendu, lorsque nous parlons d'un topos annelé $(X,A)$ sans préciser, que $A$ est commutatif.}

\subsection{}\label{higgs0-not1}
Dans cette introduction, $K$ désigne un corps de valuation discrète complet de caractéristique
$0$, à corps résiduel parfait $k$ de caractéristique $p>0$, $\co_K$ l'anneau de valuation de $K$, 
$\oK$ une clôture algébrique de $K$, $\co_\oK$ la clôture intégrale de $\co_K$ dans $\oK$, 
$\co_C$ le séparé complété $p$-adique de $\co_\oK$ et $C$ le corps des fractions de $\co_C$. 
À partir de § \ref{higgs0-rmd}, on supposera $k$ algébriquement clos. 
On pose $S=\Spec(\co_K)$, $\oS=\Spec(\co_\oK)$ et $\coS=\Spec(\co_C)$. 
On note $s$ (resp. $\eta$, resp. $\oeta$) le point fermé de $S$ (resp. générique de $S$, resp. générique de $\oS$).
Pour tout entier $n\geq 1$ et tout $S$-schéma $X$, on pose $S_n=\Spec(\co_K/p^n\co_K)$, 
\begin{equation}\label{higgs0-not1a}
X_n=X\times_SS_n,\ \ \ \oX=X\times_S\oS\ \ \ {\rm et} \ \ \ \coX=X\times_S\coS.
\end{equation}
Pour tout groupe abélien $M$,  on note $\hM$ son séparé complété $p$-adique. 

\subsection{}\label{higgs0-not15}
Soient $G$ un groupe profini, $A$ un anneau topologique muni d'une action continue de $G$ 
par des homomorphismes d'anneaux. Une {\em $A$-représentation} de $G$ est la donnée d'un $A$-module 
$M$ et d'une action $A$-semi-linéaire de $G$ sur $M$, {\em i.e.},  telle que 
pour tous $g\in G$, $a\in A$ et $m\in M$, on ait $g(am)=g(a)g(m)$.
On dit que la $A$-représentation est {\em continue} si $M$ est un $A$-module topologique 
et si l'action de $G$ sur $M$ est continue. Soient $M$, $N$ deux $A$-représentations 
(resp. deux $A$-représentations continues) de $G$. 
Un morphisme de $M$ dans $N$ est la donnée d'un morphisme $A$-linéaire et $G$-équivariant 
(resp. $A$-linéaire, continu et $G$-équivariant) de $M$ dans $N$.

\subsection{}\label{higgs0-not2}
Soit $(X,A)$ un topos annelé, $E$ un $A$-module.
Un {\em $A$-module de Higgs à coefficients dans $E$} est un couple $(M,\theta)$ formé 
d'un $A$-module $M$ et d'un morphisme $A$-linéaire $\theta\colon M\rightarrow M\otimes_AE$
tel que $\theta\wedge \theta=0$ (cf. \cite{ag1} 2.8). Suivant Simpson (\cite{simpson4} page 24),
on appelle complexe de {\em Dolbeault} de $(M,\theta)$ et 
l'on note $\mK^\bullet(M,\theta)$ le complexe de cochaînes de $A$-modules 
\begin{equation}\label{higgs0-not2a}
M\rightarrow M\otimes_AE\rightarrow M\otimes_A\wedge^2E\dots
\end{equation}
déduit de $\theta$ (cf. \cite{ag1} 2.8.2).  

\subsection{}\label{higgs0-not6}
Soient $(X,A)$ un topos annelé, $B$ une $A$-algèbre, $M$ un $B$-module, $\lambda\in \Gamma(X,A)$. 
Une {\em $\lambda$-connexion} sur $M$ relativement à l'extension $B/A$  
est la donnée d'un morphisme $A$-linéaire 
\begin{equation}
\nabla\colon M\rightarrow \Omega^1_{B/A}\otimes_BM
\end{equation}
tel que pour toutes sections locales $x$ de $B$ et $s$ de $M$, on ait 
\begin{equation}\label{higgs0-not6a}
\nabla(xs)=\lambda d(x)\otimes s+x\nabla(s).
\end{equation} 
Elle est {\em intégrable} si $\nabla \circ \nabla=0$ (cf. \cite{ag1} 2.10). 
On omettra l'extension $B/A$ de la terminologie lorsqu'il n'y a aucun risque de confusion. 

Soient $(M,\nabla)$, $(M',\nabla')$ deux modules à $\lambda$-connexions. 
Un morphisme de $(M,\nabla)$ dans $(M',\nabla')$
est la donnée d'un morphisme $B$-linéaire $u\colon M\rightarrow M'$ 
tel que  $(\id \otimes u)\circ \nabla=\nabla'\circ u$. 

Les $1$-connexions sont classiquement appelées {\em connexions}. 
Les $0$-connexions intégrables sont les $B$-champs de Higgs à coefficients dans $\Omega^1_{B/A}$.

\begin{rema}\label{higgs0-not7}
Soient $(X,A)$ un topos annelé, $B$ une $A$-algèbre, $\lambda\in \Gamma(X,A)$, $(M,\nabla)$ un module 
à $\lambda$-connexion relativement à l'extension $B/A$. Supposons qu'il existe un $A$-module $E$
et un isomorphisme $B$-linéaire $\gamma\colon E\otimes_AB\stackrel{\sim}{\rightarrow}\Omega^1_{B/A}$ 
tels que pour toute section locale $\omega$ de $E$, on ait $d(\gamma(\omega\otimes 1))=0$.
Pour que la $\lambda$-connexion $\nabla$ soit intégrable, il faut et il suffit que 
le morphisme $\theta\colon M\rightarrow E\otimes_AM$ induit par $\nabla$ et $\gamma$,
soit un $A$-champ de Higgs sur $M$ à coefficients dans $E$ (cf. \cite{ag1} 2.12). 
\end{rema}

\subsection{}\label{higgs0-not3}
Si $\cC$ est une catégorie additive, on désigne par $\cC_\mQ$ et l'on appelle 
{\em catégorie des objets de $\cC$ à isogénie près}, la catégorie ayant mêmes objets que $\cC$, 
et telle que l'ensemble des morphismes entre deux objets soit donné par 
\begin{equation}\label{higgs0-not3a}
\Hom_{\cC_\mQ}(E,F)=\Hom_{\cC}(E,F)\otimes_\mZ\mQ.
\end{equation}
La catégorie $\cC_\mQ$ n'est autre que la catégorie localisée de $\cC$ par rapport au système multiplicatif 
des {\em isogénies} de $\cC$ (cf. \cite{ag3} 6.1). On désigne par
\begin{equation}\label{higgs0-not3b}
\cC\rightarrow \cC_\mQ, \ \ \ M\mapsto M_\mQ
\end{equation}
le foncteur de localisation. 

Si $\cC$ est une catégorie abélienne, la catégorie $\cC_\mQ$ est abélienne et le foncteur de localisation 
\eqref{higgs0-not3b} est exact. En fait, $\cC_\mQ$ s'identifie canoniquement
à la catégorie quotient de $\cC$ par la sous-catégorie épaisse des objets d'exposant fini (\cite{ag3} 6.1.4).

\subsection{}\label{higgs0-not4}
Soit $(X,A)$ un topos annelé. On désigne par $\bMod(A)$ la catégorie des $A$-modules
de $X$ et par $\bMod_\mQ(A)$, au lieu de $\bMod(A)_\mQ$, la catégorie des $A$-modules à isogénie près \eqref{higgs0-not3}.
Le produit tensoriel des $A$-modules induit un bifoncteur 
\begin{equation}\label{higgs0-not4a}
\bMod_{\mQ}(A)\times \bMod_{\mQ}(A)\rightarrow \bMod_{\mQ}(A),\ \ \ (M,N)\mapsto M\otimes_{A_\mQ}N,
\end{equation} 
faisant de $\bMod_{\mQ}(A)$ une catégorie monoïdale symétrique, ayant $A_\mQ$ pour objet unité. 
Les objets de $\bMod_{\mQ}(A)$ seront aussi appelés des {\em $A_\mQ$-modules}. Cette terminologie
se justifie en considérant $A_\mQ$ comme un monoïde de $\bMod_{\mQ}(A)$.

\subsection{}\label{higgs0-imh1}
Soient $(X,A)$ un topos annelé, $E$ un $A$-module. 
On appelle {\em $A$-isogénie de Higgs à coefficients dans $E$} la donnée d'un quadruplet
\begin{equation}\label{higgs0-imh1a}
(M,N,u\colon M\rightarrow N,\theta\colon M\rightarrow N\otimes_AE)
\end{equation}
formé de deux $A$-modules $M$ et $N$ et de deux morphismes $A$-linéaires $u$ et $\theta$  
vérifiant la propriété suivante~: il existe un entier $n\not=0$ et un morphisme $A$-linéaire $v\colon N\rightarrow M$ tels que 
$v\circ u=n\cdot \id_M$, $u\circ v=n\cdot \id_N$ et que $(M,(v\otimes \id_E)\circ \theta)$ et $(N,\theta\circ v)$ 
soient des $A$-modules de Higgs à coefficients dans $E$ \eqref{higgs0-not2}. 
On notera que $u$ induit une isogénie de modules de Higgs de $(M,(v\otimes \id_E)\circ \theta)$
dans $(N,\theta\circ v)$  (\cite{ag3} 6.1), d'où la terminologie.  
Soient $(M,N,u,\theta)$, $(M',N',u',\theta')$ deux $A$-isogénies de Higgs à coefficients dans $E$. 
Un morphisme de $(M,N,u,\theta)$ dans $(M',N',u',\theta')$ est la donnée de deux morphismes $A$-linéaires
$\alpha\colon M\rightarrow M'$ et $\beta\colon N\rightarrow N'$ tels que $\beta\circ u=u'\circ \alpha$
et $(\beta\otimes \id_E)\circ \theta=\theta'\circ \alpha$. 
On désigne par $\bIH(A,E)$ la catégorie des $A$-isogénies de Higgs à coefficients dans $E$. C'est une catégorie 
additive. On note $\bIH_\mQ(A,E)$ la catégorie des objets de $\bIH(A,E)$ à isogénie près. 

\subsection{}\label{higgs0-not5}
Soient $(X,A)$ un topos annelé, $B$ une $A$-algèbre, $\lambda\in \Gamma(X,A)$. 
On appelle {\em $\lambda$-isoconnexion relativement à l'extension $B/A$}
(ou simplement {\em $\lambda$-isoconnexion} lorsqu'il n'y a aucun risque de confusion) la donnée d'un quadruplet   
\begin{equation}\label{higgs0-not5a}
(M,N,u\colon M\rightarrow N,\nabla\colon M\rightarrow \Omega^1_{B/A}\otimes_BN)
\end{equation}
où $M$ et $N$ sont des $B$-modules, $u$ est une isogénie de $B$-modules (\cite{ag3} 6.1)
et $\nabla$ est un morphisme $A$-linéaire tel que pour toutes sections locales 
$x$ de $B$ et $t$ de $M$, on ait 
\begin{equation}\label{higgs0-not5b}
\nabla(xt)=\lambda d(x) \otimes u(t)+x\nabla(t).
\end{equation} 
Pour tout morphisme $B$-linéaire $v\colon N\rightarrow M$ pour lequel il existe un entier $n$ 
tels que  $u\circ v=n\cdot \id_N$ et $v\circ u=n\cdot \id_M$, 
les couples $(M,(\id\otimes v)\circ \nabla)$ et $(N,\nabla\circ v)$ sont des modules à $(n\lambda)$-connexions \eqref{higgs0-not15},
et $u$ est un morphisme de $(M,(\id\otimes v)\circ \nabla)$ dans $(N,\nabla\circ v)$. On dit que la $\lambda$-isoconnexion 
$(M,N,u,\nabla)$ est {\em intégrable} s'il existe un morphisme $B$-linéaire $v\colon N\rightarrow M$ et un entier $n\not= 0$ 
tels que  $u\circ v=n\cdot \id_N$, $v\circ u=n\cdot \id_M$ et que les $(n\lambda)$-connexions
$(\id\otimes v)\circ \nabla$ sur $M$ et $\nabla\circ v$ sur $N$ soient intégrables. 

Soient $(M,N,u,\nabla)$, $(M',N',u',\nabla')$ deux $\lambda$-isoconnexions. 
Un morphisme de $(M,N,u,\nabla)$ dans $(M',N',u',\nabla')$ est la donnée de 
deux morphismes $B$-linéaires $\alpha\colon M\rightarrow M'$ et $\beta\colon N\rightarrow N'$
tels que  $\beta\circ u=u'\circ \alpha$ et $(\id \otimes \beta)\circ \nabla=\nabla'\circ \alpha$.

\section{Petites représentations généralisées}\label{higgs0-prg}

\subsection{}\label{higgs0-prg1}
Dans cette section, on se donne un $S$-schéma affine lisse $X=\Spec(R)$ 
tel que $X_\oeta$ soit connexe et que $X_s$ soit non-vide, un entier $d\geq 1$ et un $S$-morphisme étale 
\begin{equation}\label{higgs0-prg1a}
X\rightarrow \mG_{m,S}^d=\Spec(\co_K[T_1^{\pm 1},\dots, T_d^{\pm 1}]).
\end{equation}
C'est l'exemple type de petit schéma affine de Faltings. 
L'hypothèse de connexité de $X_\oeta$ n'est pas nécessaire mais permet de simplifier la présentation.
Le lecteur reconnaitra la nature logarithmique de la donnée \eqref{higgs0-prg1a}. 
D'ailleurs, suivant \cite{faltings3}, nous considérons 
dans ce travail une situation logarithmique lisse plus générale, qui s'avère nécessaire 
même pour définir la correspondance de Simpson $p$-adique pour une courbe propre et lisse sur $S$.
En effet, dans la seconde étape de descente, nous serons amenés à considérer des revêtements finis de 
sa fibre générique, ce qui nous ramène au cas d'un schéma semi-stable sur $S$. 
Toutefois, pour simplifier la présentation, nous nous limitons dans cette introduction au cas lisse au sens usuel
(cf. \cite{ag1} 6.2 pour le cas affine logarithmique lisse).   
On note $t_i$ l'image de $T_i$ dans $R$ $(1\leq i\leq d)$ et on pose 
\begin{equation}\label{higgs0-prg1b}
R_1=R\otimes_{\co_K}\co_{\oK}.
\end{equation}

\subsection{}\label{higgs0-prg2}
Soient $\oy$ un point géométrique de $X_\oeta$, $(V_i)_{i\in I}$ un revêtement universel de $X_\oeta$
en $\oy$. Notons $\Delta$ le groupe fondamental géométrique $\pi_1(X_\oeta,\oy)$. 
Pour tout $i\in I$, on désigne par $\oX_i=\Spec(R_i)$ la clôture intégrale de $\oX$ dans $V_i$ et on pose
\begin{equation}\label{higgs0-prg2a}
\oR=\underset{\underset{i\in I}{\longrightarrow}}{\lim}\ R_i.
\end{equation}
Dans ce contexte, les {\em représentations généralisées} de $\Delta$ sont 
les $\hoR$-représentations continues de $\Delta$ à valeurs dans des $\hoR$-modules projectifs de type fini, 
munis de leurs topologies $p$-adiques \eqref{higgs0-not15}. 
Une telle représentation $M$ est dite {\em petite} si $M$ est un $\hoR$-module libre de type fini 
ayant une base formée d'éléments $\Delta$-invariants modulo $p^{2\alpha}M$ 
pour un nombre rationnel $\alpha>\frac{1}{p-1}$.
La principale vertu des petites représentations généralisées de $\Delta$ 
est leur bon comportement par descente pour certains quotients de $\Delta$ isomorphes à $\mZ_p(1)^d$. 
Fixons un tel quotient $\Delta_\infty$ en choisissant, pour tout $1\leq i\leq d$, 
un système compatible $(t_i^{(n)})_{n\in \mN}$ de racines d'ordre $p^n$ de $t_i$ dans $\oR$.  
On définit de même la notion de $\hRun$-représentation petite de $\Delta_\infty$. Le foncteur 
\begin{equation}\label{higgs0-prg2b}
M\mapsto M\otimes_{\hRun}\hoR 
\end{equation}
de la catégorie des petites $\hRun$-représentations de $\Delta_\infty$ dans celle des 
petites $\hoR$-représentations de $\Delta$ est alors une équivalence de catégories (cf. \cite{ag1} 14.4). 
C'est une conséquence du profond théorème de presque pureté de Faltings (cf. \cite{ag1} 6.16; \cite{faltings2} § 2b).

\subsection{}\label{higgs0-prg3}
Si $(M,\varphi)$ est une petite $\hRun$-représentation de $\Delta_\infty$, 
on peut considérer le logarithme de $\varphi$ qui est un homomorphisme de $\Delta_\infty$ 
dans $\End_{\hRun}(M)$. Fixant une $\mZ_p$-base $\zeta$ de $\mZ_p(1)$, 
ce dernier s'écrit de manière unique sous la forme 
\begin{equation}\label{higgs0-prg3a}
\log(\varphi)=\sum_{i=1}^d\theta_i\otimes \chi_i\otimes \zeta^{-1},
\end{equation}
où $\zeta^{-1}$ est la base duale de $\mZ_p(-1)$, 
$\chi_i$ est le caractère de $\Delta_\infty$ à valeurs dans $\mZ_p(1)$ qui donne son action 
sur le système $(t_i^{(n)})_{n\in \mN}$ et $\theta_i$ est un endomorphisme $\hRun$-linéaire de $M$. 
Il est immédiat de voir que 
\begin{equation}\label{higgs0-prg3b}
\theta=\sum_{i=1}^d\theta_i\otimes d\log(t_i)\otimes \zeta^{-1}
\end{equation}
est un $\hRun$-champ de Higgs sur $M$ à coefficients dans $\Omega^1_{R/\co_K}\otimes_R\hR_1(-1)$
\eqref{higgs0-not2} (nous dirons  pour simplifier à coefficients dans $\Omega^1_{R/\co_K}(-1)$). 
La correspondance $(M,\varphi)\mapsto (M,\theta)$ ainsi définie est en fait une {\em équivalence de catégories}
entre la catégorie des petites $\hRun$-représentations de $\Delta_\infty$ et celle des {\em petits} $\hRun$-modules 
de Higgs à coefficients dans $\Omega^1_{R/\co_K}(-1)$ (c'est-à-dire, la catégorie des $\hRun$-modules de Higgs à 
coefficients dans $\Omega^1_{R/\co_K}(-1)$ dont le $\hRun$-module sous-jacent est libre de type fini et dont 
le champ de Higgs est nul modulo $p^{2\alpha}$ pour un nombre rationnel $\alpha>\frac{1}{p-1}$). 
Combinée avec l'énoncé précédent de descente \eqref{higgs0-prg2b}, on obtient une équivalence entre 
la catégorie des petites $\hoR$-représentations de $\Delta$ et celle des petits $\hRun$-modules 
de Higgs à coefficients dans $\Omega^1_{R/\co_K}(-1)$. L'inconvénient de cette construction est sa dépendance
en les $(t_i^{(n)})_{n\in \mN}$ $(1\leq i\leq d)$, qui exclut toute globalisation. Pour remédier à ce défaut, 
Faltings propose une autre définition, équivalente, mais dépendante d'un autre choix qui se globalise facilement.   
Notre approche, qui fait l'objet de la suite de cette introduction, est inspirée de cette construction.  

\section{Le torseur des déformations}\label{higgs0-dth}

\subsection{}\label{higgs0-dth0}
Dans cette section, on se donne un $S$-schéma affine lisse $X=\Spec(R)$ 
tel que $X_\oeta$ soit connexe, que $X_s$ soit non-vide,
et qu'il existe un entier $d\geq 1$ et un $S$-morphisme étale 
$X\rightarrow \mG_{m,S}^d$ (mais on ne fixe pas un tel morphisme). 
On se donne aussi un point géométrique $\oy$ de $X_\oeta$ et un revêtement universel $(V_i)_{i\in I}$  
de $X_\oeta$ en $\oy$, et on reprend les notations de \ref{higgs0-prg2}~: 
$\Delta=\pi_1(X_\oeta,\oy)$, $R_1=R\otimes_{\co_K}\co_{\oK}$ et $\oR$ \eqref{higgs0-prg2a}.

\subsection{}\label{higgs0-dth1}
Rappelons que Fontaine associe fonctoriellement à toute $\mZ_{(p)}$-algèbre $A$ l'anneau 
\begin{equation}\label{higgs0-dth1a}
\cR_A=\underset{\underset{x\mapsto x^p}{\longleftarrow}}{\lim}A/pA,
\end{equation} 
et un homomorphisme $\theta$ de l'anneau $\rW(\cR_A)$
des vecteurs de Witt de $\cR_A$ dans le séparé complété $p$-adique $\hA$ de $A$ (cf. \cite{ag1} 9.3). 
On pose 
\begin{equation}\label{higgs0-dth1b}
\cA_2(A)=\rW(\cR_A)/\ker(\theta)^2
\end{equation} 
et on note encore $\theta\colon \cA_2(\cR_A)\rightarrow \hA$ l'homomorphisme
induit par $\theta$. 

Dans la suite de cet article, on fixe une suite 
$(p_n)_{n\in \mN}$ d'éléments de $\co_\oK$ telle que $p_0=p$ et $p_{n+1}^p=p_n$ pour tout $n\geq 0$. 
On désigne par $\upp$ l'élément de $\cR_{\co_\oK}$ induit par la suite $(p_n)_{n\in \mN}$ et on pose  
\begin{equation}\label{higgs0-dth1d}
\xi=[\upp]-p \in \rW(\cR_{\co_\oK}),
\end{equation}
où $[\ ]$ est le représentant multiplicatif. La suite 
\begin{equation}\label{higgs0-dth1e}
0\longrightarrow \rW(\cR_{\co_\oK})\stackrel{\cdot \xi}{\longrightarrow} \rW(\cR_{\co_\oK})
\stackrel{\theta}{\longrightarrow} \co_C \longrightarrow 0
\end{equation}
est exacte (\cite{ag1} 9.5). Elle induit une suite exacte 
\begin{equation}\label{higgs0-dth1f}
0\longrightarrow \co_C\stackrel{\cdot \xi}{\longrightarrow} \cA_2(\co_\oK)
\stackrel{\theta}{\longrightarrow} \co_C \longrightarrow 0,
\end{equation}
où on a encore noté $\cdot \xi$ le morphisme déduit de la multiplication par $\xi$ dans $\cA_2(\co_\oK)$. 
L'idéal $\ker(\theta)$ de $\cA_2(\co_\oK)$ est de carré nul. 
C'est un $\co_C$-module libre de base $\xi$. Il sera noté $\xi\co_C$. 
On observera que contrairement à $\xi$, ce module ne dépend pas du choix de la suite $(p_n)_{n\in \mN}$. 
On note $\xi^{-1}\co_C$ le $\co_C$-module dual de $\xi\co_C$. 
Pour tout $\co_C$-module $M$, on désigne les $\co_C$-modules $M\otimes_{\co_C}(\xi \co_C)$ 
et $M\otimes_{\co_C}(\xi^{-1} \co_C)$ simplement par $\xi M$ et $\xi^{-1} M$, respectivement. 

De même, on a une suite exacte (\cite{ag1} (9.11.2))
\begin{equation}\label{higgs0-dth1c}
0\longrightarrow \hoR\stackrel{\cdot \xi}{\longrightarrow} \cA_2(\oR)\stackrel{\theta}{\longrightarrow} \hoR\longrightarrow 0.
\end{equation}
L'idéal $\ker(\theta)$ de $\cA_2(\oR)$ est de carré nul. C'est un $\hoR$-module libre de base $\xi$, 
canoniquement isomorphe à $\xi\hoR$. Le groupe $\Delta$ agit par fonctorialité sur $\cA_2(\oR)$. 

Posons $\cA_2(\oS)=\Spec(\cA_2(\co_\oK))$, $Y=\Spec(\oR)$, $\hY=\Spec(\hoR)$ et $\cA_2(Y)=\Spec(\cA_2(\oR))$.

\subsection{}\label{higgs0-dth11}
On se donne dans la suite une $\cA_2(\oS)$-déformation lisse $\tX$ de 
$\coX$, c'est-à-dire, un $\cA_2(\oS)$-schéma lisse $\tX$ qui s'insère dans un diagramme cartésien 
\begin{equation}\label{higgs0-dth11a}
\xymatrix{
{\coX}\ar[r]\ar[d]&{\tX}\ar[d]\\
{\coS}\ar[r]&{\cA_2(\oS)}}
\end{equation}
Cette donnée supplémentaire remplace la donnée d'un $S$-morphisme étale $X\rightarrow \mG_{m,S}^d$; 
d'ailleurs, un tel morphisme fournit une déformation. 

On pose 
\begin{equation}
\rT=\Hom_{\hoR}(\Omega^1_{R/\co_K}\otimes_R\hoR,\xi\hoR).
\end{equation} 
On identifie le $\hoR$-module dual à $\xi^{-1}\Omega^1_{R/\co_K}\otimes_R\hoR$ \eqref{higgs0-dth1}
et on note $\bT$ le $\hY$-fibré vectoriel associé, autrement dit,  
\begin{equation}\label{higgs0-tor2c}
\bT=\Spec(\Sym_{\hoR}(\xi^{-1}\Omega^1_{R/\co_K}\otimes_R\hoR)).
\end{equation}
Soient $U$ un ouvert de $\hY$, $\tU$ l'ouvert de $\cA_2(Y)$ défini par $U$. 
On désigne par $\cL(U)$ l'ensemble des morphismes représentés par des flèches pointillées qui complètent  le diagramme 
\begin{equation}\label{higgs0-dth11b}
\xymatrix{
{U}\ar[r]\ar[d]&{\tU}\ar@{.>}[d]\ar@/^2pc/[dd]\\
{\coX}\ar[r]\ar[d]&{\tX}\ar[d]\\
{\coS}\ar[r]&{\cA_2(\oS)}}
\end{equation}
de façon à le laisser commutatif. Le foncteur $U\mapsto \cL(U)$ est un $T$-torseur pour la topologie de Zariski de $\hY$. 
On désigne par $\cF$ le $\hoR$-module des fonctions affines sur $\cL$ (cf. \cite{ag1} 4.9). 
Celui-ci s'insère dans une suite exacte canonique (\cite{ag1} (4.9.1))
\begin{equation}\label{higgs0-dth11d}
0\rightarrow \hoR\rightarrow \cF\rightarrow \xi^{-1}\Omega^1_{R/\co_K} \otimes_R \hoR\rightarrow 0.
\end{equation} 
Cette suite induit pour tout entier $n\geq 1$ une suite exacte 
\begin{equation}\label{higgs0-dth11e}
0\rightarrow \Sym^{n-1}_{\hoR}(\cF)\rightarrow \Sym^{n}_{\hoR}(\cF)\rightarrow \Sym^n_{\hoR}(\xi^{-1}\Omega^1_{R/\co_K}
\otimes_R\hoR)\rightarrow 0.
\end{equation}
Les $\hoR$-modules $(\Sym^{n}_{\hoR}(\cF))_{n\in \mN}$ forment donc un système inductif filtrant 
dont la limite inductive 
\begin{equation}\label{higgs0-dth11f}
\cC=\underset{\underset{n\geq 0}{\longrightarrow}}\lim\ \Sym^n_{\hoR}(\cF)
\end{equation}
est naturellement munie d'une structure de $\hoR$-algèbre. D'après (\cite{ag1} 4.10), le $\hY$-schéma 
\begin{equation}\label{higgs0-dth11c}
\bL=\Spec(\cC)
\end{equation}
est naturellement un $\bT$-fibré principal homogène sur $\hY$ qui représente canoniquement $\cL$. 

L'action naturelle de $\Delta$ sur le schéma $\cA_2(Y)$ induit  
une action $\hoR$-semi-linéaire de $\Delta$ sur $\cF$, telle que les morphismes de la suite \eqref{higgs0-dth11d} soient 
$\Delta$-équivariants. On en déduit une action de $\Delta$ sur 
$\cC$ par des automorphismes d'anneaux, compatible avec son action sur $\hoR$, 
que l'on appelle {\em action canonique}. Ces actions sont continues pour les topologies $p$-adiques (\cite{ag1} 12.4).
La $\hoR$-algèbre $\cC$, munie de l'action canonique de $\Delta$, 
est appelée {\em l'algèbre de Higgs-Tate} associée à $\tX$. 
La $\hoR$-représentation $\cF$ de $\Delta$ est appelée 
l'{\em extension de Higgs-Tate} associée à $\tX$.

\subsection{}\label{higgs0-dth2}
Soit $(M,\theta)$ un {\em petit} $\hRun$-module de Higgs à coefficients dans $\xi^{-1}\Omega^1_{R/\co_K}$
(c'est-à-dire, un $\hRun$-module de Higgs à coefficients dans $\xi^{-1}\Omega^1_{R/\co_K}\otimes_R\hRun$ 
dont le $\hRun$-module sous-jacent est libre de type fini et dont 
le champ de Higgs est nul modulo $p^{\alpha}$ pour un nombre rationnel $\alpha>\frac{1}{p-1}$)
et soit $\psi\in \cL(\hY)$.
Pour tout $\sigma\in \Delta$, on note ${^\sigma \psi}$ la section de $\cL(\hY)$ définie par le diagramme commutatif
\begin{equation}
\xymatrix{
{\bL}\ar[r]^\sigma\ar[d]&{\bL}\ar[d]\\
{\hY}\ar[r]^\sigma\ar@/^1pc/[u]^{\psi}&{\hY}\ar@/_1pc/[u]_{{^\sigma \psi}}}
\end{equation}
La différence $D_\sigma=\psi-{^\sigma \psi}$
est un élément de $\Hom_{\hoR}(\xi^{-1}\Omega^1_{R/\co_K}\otimes_R\hoR,\hoR)$. 
L'endomorphisme $\exp((D_\sigma\otimes \id_M)\circ \theta)$ de $M\otimes_{\hRun}\hoR$ est bien défini, 
compte tenu de la petitesse de $\theta$. 
On obtient alors une petite $\hoR$-représentation de $\Delta$ sur $M\otimes_{\hRun}\hoR$. 
La correspondance ainsi définie est en fait une équivalence de catégories de la catégorie des petits 
$\hRun$-modules de Higgs à coefficients 
dans $\xi^{-1}\Omega^1_{R/\co_K}$ dans celle des petites $\hoR$-représentations de $\Delta$. 
C'est essentiellement un quasi-inverse de l'équivalence de catégories définie dans \ref{higgs0-prg3}.

Pour éviter le choix d'une section $\psi$ de $\cL(\hY)$, on peut effectuer le changement de base de $\hoR$
à $\cC$ et utiliser le plongement diagonal de $\bL$. La construction précédente dans ce cadre 
peut alors s'interpréter suivant le schéma classique des 
correspondances introduites par Fontaine (ou encore plus classique de la correspondance de Riemann-Hilbert 
analytique complexe) en prenant pour anneau de périodes faisant le pont entre les représentations généralisées
et les modules de Higgs un complété $p$-adique faible $\cC^\dagger$ de $\cC$
(la complétion est rendue nécessaire par l'exponentielle). 
A cet anneau est naturellement associée une notion d'admissibilité~; 
c'est la notion de {\em représentation généralisée de Dolbeault}.
Avant de développer cette approche, nous dirons quelques mots sur l'anneau $\cC$ qui peut lui-même servir d'anneau de périodes 
entre représentations généralisées et modules de Higgs.  
En fait, $\cC$ est un {\em modèle entier de l'anneau de Hyodo} (cf. \eqref{higgs0-dth4c} et \cite{ag1} 15.6), ce qui explique 
le lien entre notre approche et celle de Hyodo.

\subsection{}\label{higgs0-dth3}
On rappelle que Faltings a défini une extension canonique de $\hoR$-représentations de $\pi_1(X,\oy)$
\begin{equation}\label{higgs0-dth3d}
0\rightarrow \rho^{-1}\hoR\rightarrow \cE\rightarrow \Omega^1_{R/\co_K}\otimes_R\hoR(-1)\rightarrow 0,
\end{equation}
où $\rho$ est un élément de $\co_\oK$ de valuation $\geq \frac{1}{p-1}$,
qui joue un rôle important dans son approche de la théorie de Hodge $p$-adique  (cf. \cite{ag1} 7.22). 
Nous montrons dans (\cite{ag1} 10.19) qu'il existe un morphisme $\hoR$-linéaire et $\Delta$-équivariant 
\begin{equation}\label{higgs0-dth3e}
p^{-\frac{1}{p-1}}\cF\rightarrow \cE
\end{equation}
qui s'insère dans un diagramme commutatif 
\begin{equation}\label{higgs0-dth3f}
\xymatrix{
0\ar[r]&{p^{-\frac{1}{p-1}}\hoR}\ar[r]\ar@{^(->}[d]&{p^{-\frac{1}{p-1}}\cF}\ar[r]\ar[d]&
{p^{-\frac{1}{p-1}}\xi^{-1}\Omega^1_{R/\co_K} \otimes_R \hoR}\ar[r]\ar[d]^{-c}&0\\
0\ar[r]&{\rho^{-1}\hoR}\ar[r]&{\cE}\ar[r]&{\Omega^1_{R/\co_K} \otimes_R \hoR(-1)}\ar[r] & 0}
\end{equation}
où $c$ est l'isomorphisme induit par un isomorphisme canonique 
$\hoR(1)\stackrel{\sim}{\rightarrow}p^{\frac{1}{p-1}}\xi\hoR$
(\cite{ag1} 9.18). Le morphisme \eqref{higgs0-dth3e} est canonique si l'on prend pour $\tX$ 
la déformation induite par un $S$-morphisme étale $X\rightarrow \mG_{m,S}^d$. 
Il est important de noter que dans le cadre logarithmique 
qui sera considéré dans ce travail, l'extension de Faltings change légèrement de forme puisque le facteur 
$\rho^{-1}\hoR$ est remplacé par $(\pi\rho)^{-1}\hoR$, où $\pi$ est une uniformisante de $K$.

\subsection{}\label{higgs0-dth4}
Partant de l'extension $\cE$ de Faltings \eqref{higgs0-dth3d}, Hyodo \cite{hyodo1} définit une $\hoR$-algèbre $\cC_\HT$ 
par une limite inductive analogue à \eqref{higgs0-dth11f}. On notera que $p$ étant inversible dans $\cC_\HT$, il revient au 
même de partir de $\cE\otimes_{\mZ_p}\mQ_p$, ce qui correspond à la définition originelle de Hyodo. 
Le morphisme \eqref{higgs0-dth3e} induit donc un isomorphisme $\Delta$-équivariant de $\hoR$-algèbres 
\begin{equation}\label{higgs0-dth4c}
\cC[\frac 1 p]\stackrel{\sim}{\rightarrow} \cC_\HT.
\end{equation}
Pour toute $\mQ_p$-représentation continue $V$ de $\Gamma=\pi_1(X,\oy)$ et tout entier $i$, 
Hyodo définit le $\hR[\frac 1 p]$-module $\rD^i(V)$ par  
\begin{equation}\label{higgs0-dth4a}
\rD^i(V)=(V\otimes_{\mQ_p}\cC_\HT(i))^\Gamma.
\end{equation}
La représentation $V$ est dite de {\em Hodge-Tate} si elle satisfait aux conditions suivantes~:  
\begin{itemize}
\item[{\rm (i)}] $V$ est un $\mQ_p$-espace vectoriel de dimension finie, muni de la topologie $p$-adique. 
\item[{\rm (ii)}] Le morphisme canonique 
\begin{equation}\label{higgs0-dth4b}
\oplus_{i\in \mZ}\rD^i(V)\otimes_{\hR[\frac 1 p]}\cC_\HT(-i)\rightarrow V\otimes_{\mQ_p}\cC_\HT
\end{equation}
est un isomorphisme. 
\end{itemize}

\subsection{}\label{higgs0-dolbeault1}
Pour tout nombre rationnel $r\geq 0$, on note $\cF^{(r)}$ la $\hoR$-représentation de $\Delta$ 
déduite de $\cF$ par image inverse par la multiplication par $p^r$ sur 
$\xi^{-1}\Omega^1_{R/\co_K}\otimes_R\hoR$, de sorte qu'on a une suite exacte
\begin{equation}\label{higgs0-dolbeault1a}
0\rightarrow \hoR\rightarrow \cF^{(r)}\rightarrow 
\xi^{-1}\Omega^1_{R/\co_K}\otimes_R\hoR\rightarrow 0.
\end{equation}
Cette suite induit pour tout entier $n\geq 1$, une suite exacte 
\begin{equation}\label{higgs0-dolbeault1b}
0\rightarrow \Sym^{n-1}_{\hoR}(\cF^{(r)})\rightarrow \Sym^{n}_{\hoR}(\cF^{(r)})\rightarrow 
\Sym^n_{\hoR}(\xi^{-1}\Omega^1_{R/\co_K}
\otimes_R\hoR)\rightarrow 0.
\end{equation}
Les $\hoR$-modules $(\Sym^{n}_{\hoR}(\cF^{(r)}))_{n\in \mN}$ forment donc un système inductif filtrant, 
dont la limite inductive 
\begin{equation}\label{higgs0-dolbeault1c}
\cC^{(r)}=\underset{\underset{n\geq 0}{\longrightarrow}}\lim\ \rS^n_{\hoR}(\cF^{(r)})
\end{equation}
est naturellement munie d'une structure de $\hoR$-algèbre. 
L'action de $\Delta$ sur $\cF^{(r)}$ induit une action sur $\cC^{(r)}$ 
par des automorphismes d'anneaux, compatible avec son action sur $\hoR$, que l'on appelle {\em action canonique}.
La $\hoR$-algèbre $\cC^{(r)}$ munie de cette action
est appelée {\em l'algèbre de Higgs-Tate d'épaisseur $r$} associée à $\tX$. 
On note $\hcC^{(r)}$ le séparé complété $p$-adique de $\cC^{(r)}$ que l'on suppose toujours muni de 
la topologie $p$-adique. 

Pour tous nombres rationnels $r\geq r'\geq 0$, on a un $\hoR$-homomorphisme canonique injectif et 
$\Delta$-équivariant $\alpha^{r,r'}\colon \cC^{(r')}\rightarrow \cC^{(r)}$. On vérifie aussitôt que l'homomorphisme 
induit $\halpha^{r,r'}\colon\hcC^{(r')}\rightarrow \hcC^{(r)}$ est injectif. On pose 
\begin{equation}\label{higgs0-dolbeault1f}
\cC^\dagger=\underset{\underset{r\in \mQ_{>0}}{\longrightarrow}}{\lim} \hcC^{(r)},
\end{equation}
que l'on identifie à une sous-$\hoR$-algèbre de $\hcC=\hcC^{(0)}$. 
Le groupe $\Delta$ agit naturellement sur $\cC^\dagger$ par des automorphismes d'anneaux, 
d'une façon compatible avec ses actions sur $\hoR$ et sur $\hcC$.

On désigne par  
\begin{equation}\label{higgs0-dolbeault1g}
d_{\cC^{(r)}}\colon \cC^{(r)}\rightarrow \xi^{-1}\Omega^1_{R/\co_K}\otimes_R\cC^{(r)}
\end{equation}
la $\hoR$-dérivation universelle de $\cC^{(r)}$ et par 
\begin{equation}\label{higgs0-dolbeault1h}
d_{\hcC^{(r)}}\colon \hcC^{(r)}\rightarrow \xi^{-1}\Omega^1_{R/\co_K}\otimes_R\hcC^{(r)}
\end{equation}
son prolongement aux complétés (on notera que le $R$-module $\Omega^1_{R/\co_K}$ est libre de type fini). 
Les dérivations $d_{\cC^{(r)}}$ et $d_{\hcC^{(r)}}$ sont $\Delta$-équivariantes. 
Elles sont également des $\hoR$-champs de Higgs à coefficients dans $\xi^{-1}\Omega^1_{R/\co_K}$ 
puisque $\xi^{-1}\Omega^1_{R/\co_K}\otimes_R \hoR=d_{\cC^{(r)}}(\cF^{(r)})\subset d_{\cC^{(r)}}(\cC^{(r)})$ (cf. \ref{higgs0-not7}).

Pour tous nombres rationnels $r\geq r'\geq 0$, on a 
\begin{equation}\label{higgs0-dolbeault1j}
p^{r'}(\id \times \alpha^{r,r'}) \circ d_{\cC^{(r')}}=p^rd_{\cC^{(r)}}\circ \alpha^{r,r'}.
\end{equation}
Les dérivations $p^rd_{\hcC^{(r)}}$ induisent donc une $\hoR$-dérivation 
\begin{equation}\label{higgs0-dolbeault1l}
d_{\cC^\dagger}\colon \cC^\dagger\rightarrow \xi^{-1}\Omega^1_{R/\co_K} \otimes_R\cC^\dagger,
\end{equation}
qui n'est autre que la restriction de $d_{\hcC}$ à $\cC^\dagger$. 

\subsection{}\label{higgs0-dth5}
Pour toute $\hoR$-représentation $M$ de $\Delta$, 
on note $\mH(M)$ le $\hRun$-module de Higgs à coefficients dans $\xi^{-1}\Omega^1_{R/\co_K}$ défini par 
\begin{equation}\label{higgs0-dth5a}
\mH(M)=(M\otimes_{\hoR}\cC^\dagger)^\Delta
\end{equation}
et le champ de Higgs induit par $d_{\cC^\dagger}$. 
Pour tout $\hRun$-module de Higgs $(N,\theta)$ à coefficients dans $\xi^{-1}\Omega^1_{R/\co_K}$, 
on désigne par $\mV(N)$ la $\hoR$-représentation de $\Delta$ définie par 
\begin{equation}\label{higgs0-dth5b}
\mV(N)=(N\otimes_{\hRun}\cC^\dagger)^{\theta_\tot=0},
\end{equation}
où $\theta_\tot=\theta\otimes \id+\id\otimes d_{\cC^\dagger}$, et l'action de $\Delta$ induite par son action 
canonique sur $\cC^\dagger$. Pour pouvoir exploiter pleinement ces foncteurs, 
nous établissons des résultats d'acyclicité de $\cC^\dagger\otimes_{\mZ_p}\mQ_p$ 
pour la cohomologie de Dolbeault (\cite{ag1} 12.3) et pour la cohomologie continue de $\Delta$ (\cite{ag1} 12.5),
généralisant légèrement des résultats antérieurs de Tsuji \cite{tsuji3}.  

Une $\hoR$-représentation continue $M$ de $\Delta$ 
est dite {\em de Dolbeault} si elle satisfait aux conditions suivantes (cf. \cite{ag1} 12.11)~:
\begin{itemize}
\item[{\rm (i)}] $M$ est un $\hoR$-module projectif de type fini, muni de la topologie $p$-adique~;
\item[{\rm (ii)}] $\mH(M)$ est un $\hRun$-module projectif de type fini~; 
\item[{\rm (iii)}] le morphisme $\cC^\dagger$-linéaire canonique 
\begin{equation}
\mH(M) \otimes_{\hRun}\cC^\dagger\rightarrow  M\otimes_{\hoR}\cC^\dagger
\end{equation}
est un isomorphisme.
\end{itemize}

Un $\hRun$-module de Higgs $(N,\theta)$ à coefficients dans $\xi^{-1}\Omega^1_{R/\co_K}$ est dit {\em soluble} 
s'il satisfait aux conditions suivantes (cf. \cite{ag1} 12.12)~:
\begin{itemize}
\item[{\rm (i)}] $N$ est un $\hRun$-module projectif de type fini~;
\item[{\rm (ii)}] $\mV(N)$ est un $\hoR$-module projectif de type fini~;
\item[{\rm (iii)}] le morphisme $\cC^\dagger$-linéaire canonique 
\begin{equation}
\mV(N) \otimes_{\hoR}\cC^\dagger\rightarrow  N\otimes_{\hRun}\cC^\dagger
\end{equation}
est un isomorphisme.
\end{itemize}
Il est immédiat de voir que les foncteurs $\mV$ et $\mH$ induisent des équivalences de 
catégories quasi-inverses l'une de l'autre,
entre la catégorie des $\hoR$-représentations de Dolbeault de $\Delta$ et celle des $\hRun$-modules de Higgs
solubles à coefficients dans $\xi^{-1}\Omega^1_{R/\co_K}$ (\cite{ag1} 12.15). 

Nous montrons que les petites $\hoR$-représentations de $\Delta$ sont de Dolbeault (\cite{ag1} 14.6), 
que les petits $\hRun$-modules de Higgs sont solubles (\cite{ag1} 13.20), et que $\mV$ et $\mH$ 
induisent des équivalences de catégories quasi-inverses l'une de l'autre entre les catégories de ces 
objets (\cite{ag1} 14.7). On retrouve en fait la correspondance définie dans \ref{higgs0-prg3},
à une renormalisation près (cf. \cite{ag1} 13.18). 

\subsection{}\label{higgs0-dth6}
On définit les notions de $\hoR[\frac 1 p]$-représentation de {\em Dolbeault} de $\Delta$ 
et de $\hRun[\frac 1 p]$-module de Higgs {\em soluble} à coefficients dans $\xi^{-1}\Omega^1_{R/\co_K}$ 
en calquant les définitions dans les cas entiers (cf. \cite{ag1} 12.16 et 12.18). 
Nous montrons que les foncteurs $\mV$ et $\mH$ induisent des équivalences de 
catégories quasi-inverses l'une de l'autre,
entre la catégorie des $\hoR[\frac 1 p]$-représentations de Dolbeault de $\Delta$ et celle des 
$\hRun[\frac 1 p]$-modules de Higgs solubles à coefficients dans $\xi^{-1}\Omega^1_{R/\co_K}$ (\cite{ag1} 12.24). 
Ce résultat est légèrement plus délicat que son analogue entier \eqref{higgs0-dth5}. 

Contrairement au cas entier, les conditions d'admissibilité rationnelles peuvent s'interpréter  en termes de conditions de divisibilité.
Plus précisément, on dit qu'une $\hoR[\frac 1 p]$-représentation continue $M$ de $\Delta$ est {\em petite} si
les conditions suivantes sont remplies~: 
\begin{itemize}
\item[(i)] $M$ est un $\hoR[\frac 1 p]$-module projectif de type fini, muni de la topologie $p$-adique (\cite{ag1} 2.2)~;
\item[(ii)] il existe un nombre rationnel $\alpha>\frac{2}{p-1}$
et un sous-$\hoR$-module de type fini $M^\circ$ de $M$, stable par $\Delta$, 
engendré par un nombre fini d'éléments $\Delta$-invariants modulo $p^{\alpha}M^\circ$,
et qui engendre $M$ sur $\hoR[\frac 1 p]$. 
\end{itemize}

On dit qu'un $\hRun[\frac 1 p]$-module de Higgs $(N,\theta)$ à coefficients dans $\xi^{-1}\Omega^1_{R/\co_K}$ 
est {\em petit} si les conditions suivantes sont remplies~: 
\begin{itemize}
\item[(i)] $N$ est un $\hRun[\frac 1 p]$-module projectif de type fini~; 
\item[(ii)] il existe un nombre rationnel $\beta>\frac{1}{p-1}$ et 
un sous-$\hRun$-module de type fini $N^\circ$ de $N$, qui l'engendre sur $\hRun[\frac 1 p]$, tels que l'on ait 
\begin{equation}
\theta(N^\circ)\subset p^{\beta}\xi^{-1}N^\circ \otimes_R\Omega^1_{R/\co_K}.
\end{equation}
\end{itemize}

\begin{prop}[\cite{ag1} 13.25]
Pour qu'un $\hRun[\frac 1 p]$-module de Higgs à coefficients dans $\xi^{-1}\Omega^1_{R/\co_K}$ 
soit soluble, il faut et il suffit qu'il soit petit.
\end{prop}

\begin{prop}[\cite{ag1} 13.26]
Toute $\hoR[\frac 1 p]$-représentation de Dolbeault de $\Delta$ est petite. 
\end{prop}

Nous montrons que l'implication inverse est équivalente à une propriété de descente pour les petites 
$\hoR[\frac 1 p]$-représentations de $\Delta$ (\cite{ag1} 14.8).

\begin{prop}[\cite{ag1} 12.26]\label{higgs0-dth9}
Soient $M$ une $\hoR[\frac 1 p]$-représentation de Dolbeault de $\Delta$,  
$(\mH(M),\theta)$ le $\hRun[\frac 1 p]$-module de Higgs à coefficients 
dans $\xi^{-1}\Omega^1_{R/\co_K}$ associé. 
On a alors un isomorphisme canonique fonctoriel dans $\bD^+(\bMod(\hRun[\frac 1 p]))$ 
\begin{equation}
\rC_\cont^\bullet(\Delta, M)\stackrel{\sim}{\rightarrow} \mK^\bullet(\mH(M),\theta),
\end{equation}
où $\rC_\cont^\bullet(\Delta, M)$ est le complexe de cochaînes continues de $\Delta$ dans $M$
et $\mK^\bullet(\mH(M),\theta)$ est le complexe de Dolbeault  \eqref{higgs0-not2}. 
\end{prop}

Cet énoncé a été indépendamment démontré par Faltings pour les petites représentations (\cite{faltings3} §3) et 
par Tsuji (\cite{tsuji3} 5.3.2).

\subsection{}\label{higgs0-dth7}
Il résulte de \eqref{higgs0-dth4c} que si $V$ est une $\mQ_p$-représentation 
de Hodge-Tate de $\Gamma$, alors $V\otimes_{\mZ_p}\hoR$ est une $\hoR[\frac 1 p]$-représentation de Dolbeault 
de $\Delta$; on a un isomorphisme $\hRun$-linéaire fonctoriel 
\begin{equation}\label{higgs0-dth7a}
\mH(V\otimes_{\mZ_p}\hoR)\stackrel{\sim}{\rightarrow}\oplus_{i\in \mZ}\rD^i(V)\otimes_{\hR}\hRun(-1);
\end{equation}
et le champ de Higgs sur $\mH(V\otimes_{\mZ_p}\hoR)$ est induit par les morphismes $\hR$-linéaires 
\begin{equation}\label{higgs0-dth7b}
\rD^i(V)\rightarrow \rD^{i-1}(V)\otimes_R\Omega^1_{R/\co_K}
\end{equation}
déduits de la dérivation universelle de $\cC_\HT$ sur $\hoR[\frac 1 p]$ (cf. \cite{ag1} 15.7). 
De plus, l'isomorphisme \eqref{higgs0-dth7a} est canonique si l'on prend pour $\tX$ 
la déformation induite par un $S$-morphisme étale $X\rightarrow \mG_{m,S}^d$. 

\subsection{}\label{higgs0-dth8}
Hyodo (\cite{hyodo1} 3.6) a montré que si $f\colon Y\rightarrow X$
est un morphisme propre et lisse, pour tout entier $m\geq 0$, le faisceau $\rR^mf_{\eta*}(\mQ_p)$ 
est de Hodge-Tate de poids compris entre $0$ et $m$; pour tout $0\leq i\leq m$, on a un isomorphisme 
canonique 
\begin{equation}\label{higgs0-dth7c}
\rD^i(\rR^mf_{\eta*}(\mQ_p))\stackrel{\sim}{\rightarrow}(\rR^{m-i}f_{\eta *}(\Omega^i_{Y/X}))\otimes_{R}\hR;
\end{equation}
et le morphisme \eqref{higgs0-dth7b} est induit par la classe de Kodaira-Spencer de $f$. 
Il s'ensuit que le fibré de Higgs associé à $\rR^{m}f_{\eta *}(\mQ_p)$ est égal au fibré vectoriel
\begin{equation}\label{higgs0-dth7d}
\oplus_{0\leq i\leq m} \rR^{m-i}f_{\eta *}(\Omega^{i}_{Y/X}),
\end{equation}
muni du champ de Higgs $\theta$ défini par la classe de Kodaira-Spencer de $f$.

\section{Topos annelé de Faltings}\label{higgs0-rmd}

\subsection{}\label{higgs0-rmd0}
Nous aborderons dans \cite{ag3} les aspects globaux de la théorie dans un cadre logarithmique.  
Toutefois, pour maintenir une présentation simplifiée, nous nous limitons ici encore au cas 
lisse au sens usuel  (cf. \cite{ag3} 4.7 pour le cas logarithmique lisse). 
Dans la suite de cette introduction, on suppose $k$ {\em algébriquement clos} 
et on désigne par $X$ un $S$-schéma lisse de type fini. 
\`A partir de \ref{higgs0-rmd4}, on supposera de plus qu'il existe une $\cA_2(\oS)$-déformation lisse $\tX$ de $\coX$
que l'on fixera.

\subsection{}\label{higgs0-rmd1}
La première difficulté pour recoller la 
construction locale décrite dans § \ref{higgs0-dth} est de faisceautiser la notion de représentation généralisée.
Nous utilisons pour ce faire le {\em topos de Faltings}, une variante fibrée de la notion de 
{\em topos co-évanescent} de Deligne que nous développons dans \cite{ag2}.  
On désigne par $E$ la catégorie des morphismes de schémas $V\rightarrow U$ au-dessus du morphisme 
canonique $X_\oeta\rightarrow X$, {\em i.e.}, les diagrammes commutatifs 
\begin{equation}\label{higgs0-rmd1a}
\xymatrix{
V\ar[r]\ar[d]&U\ar[d]\\
{X_\oeta}\ar[r]&X}
\end{equation}
tels que le morphisme $U\rightarrow X$ soit étale et que le morphisme $V\rightarrow U_\oeta$ soit fini étale.
Elle est fibrée au-dessus de la catégorie $\Et_{/X}$ des $X$-schémas étales, par le foncteur 
\begin{equation}\label{higgs0-rmd1b}
\pi\colon E\rightarrow \Et_{/X}, \ \ \ (V\rightarrow U)\mapsto U.
\end{equation}
La fibre de $\pi$ au-dessus d'un $X$-schéma étale $U$ est la catégorie $\Et_{\rf/U_\oeta}$ 
des schémas finis étales au-dessus
de $U_\oeta$, que l'on équipe de la topologie étale. On désigne par $U_{\oeta,\fet}$ le topos des faisceaux
d'ensembles sur $\Et_{\rf/U_\oeta}$ (cf. \cite{ag2} 9.2). 
Si $U_\oeta$ est connexe et si $\oy$ est un point géométrique de $U_\oeta$,
notant $\bB_{\pi_1(U_\oeta,\oy)}$ le topos classifiant du groupe fondamental $\pi_1(U_\oeta,\oy)$, on a une équivalence
canonique de catégories (\cite{ag2} (9.8.4))
\begin{equation}\label{higgs0-rmd1c}
\nu_\oy\colon U_{\oeta,\fet}\stackrel{\sim}{\rightarrow}\bB_{\pi_1(U_\oeta,\oy)}.
\end{equation}

Nous équipons $E$ de la {\em topologie co-évanescente} engendrée 
par les recouvrements $\{(V_i\rightarrow U_i)\rightarrow (V\rightarrow U)\}_{i\in I}$
des deux types suivants~: 
\begin{itemize}
\item[(v)] $U_i=U$ pour tout $i\in I$, et $(V_i\rightarrow V)_{i\in I}$ est un recouvrement. 
\item[(c)] $(U_i\rightarrow U)_{i\in I}$ est un recouvrement et $V_i=U_i\times_UV$ pour tout $i\in I$. 
\end{itemize}
Le site co-évanescent $E$ ainsi défini est aussi appelé {\em site de Faltings} de $X$.  
On désigne par $\tE$ et l'on appelle {\em topos de Faltings} de $X$, le topos des faisceaux d'ensembles sur $E$. 
On renvoie à \cite{ag2} pour une étude détaillée de ce topos. 
Signalons ici une description commode et simple de $\tE$.

\begin{prop}[\cite{ag2} 5.10] 
La donnée d'un faisceau $F$ sur $E$ est équivalente à la donnée pour tout objet $U$ de $\Et_{/X}$ d'un faisceau $F_U$ 
de $U_{\oeta,\fet}$ et pour tout morphisme $f\colon U'\rightarrow U$ de $\Et_{/X}$ d'un morphisme 
$F_U\rightarrow f_{\fet*}(F_{U'})$, ces morphismes étant soumis à des relations de compatibilité, 
tels que pour toute famille couvrante $(f_n\colon U_n\rightarrow U)_{n\in \Sigma}$ de $\Et_{/X}$, 
si pour tout $(m,n)\in \Sigma^2$, on pose $U_{mn}=U_m\times_UU_n$ et on note $f_{mn}\colon U_{mn}\rightarrow U$
le morphisme canonique, la suite de morphismes de faisceaux de $U_{\oeta,\fet}$
\begin{equation}\label{higgs0-rmd1d}
F_U\rightarrow \prod_{n\in \Sigma}(f_{n,\oeta})_{\fet *}(F_{U_n})\rightrightarrows 
\prod_{(m,n)\in \Sigma^2} (f_{mn,\oeta})_{\fet*}(F_{U_{mn}})
\end{equation}
soit exacte. 
\end{prop}

On identifiera dans la suite tout faisceau $F$ sur $E$ au foncteur $\{U\mapsto F_U\}$ associé, 
le faisceau $F_U$ étant la restriction de $F$ à la fibre $\Et_{\rf/U_\oeta}$ de $\pi$ au-dessus de $U$.  

\subsection{}\label{higgs0-rmd7}
Le foncteur d'injection canonique $\Et_{\rf/X_\oeta}\rightarrow E$ 
est continu et exact à gauche (\cite{ag2} 5.32). Il définit donc un morphisme de topos
\begin{equation}\label{higgs0-rmd7a}
\beta\colon \tE\rightarrow Y_\fet.
\end{equation}
De même, le foncteur 
\begin{equation}\label{higgs0-rmd7b}
\sigma^+\colon \Et_{/X}\rightarrow E, \ \ \ U\mapsto (U_\oeta\rightarrow U)
\end{equation}
est continu et exact à gauche (\cite{ag2} 5.32). Il définit donc un morphisme de topos
\begin{equation}\label{higgs0-rmd7c}
\sigma\colon \tE \rightarrow X_\et.
\end{equation}

\subsection{}\label{higgs0-rmd10}
Soient $\ox$ un point géométrique de $X$, $X'$ le localisé strict de $X$ en $\ox$. 
On désigne par $E'$ le site de Faltings associé à $X'$, 
par $\tE'$ le topos des faisceaux d'ensembles sur $E'$ et par  
\begin{equation}\label{higgs0-rmd10a}
\beta'\colon \tE'\rightarrow X'_{\oeta,\fet}
\end{equation}
le morphisme canonique \eqref{higgs0-rmd7a}. On démontre dans (\cite{ag2} 10.27) que le foncteur $\beta'_*$ est exact. 
Cette propriété est cruciale pour l'étude des principaux faisceaux du topos de Faltings 
considérés dans ce travail. 

Le morphisme canonique $X'\rightarrow X$ induit par fonctorialité un morphisme de topos (\cite{ag2} 10.12)
\begin{equation}\label{higgs0-rmd10b}
\Phi\colon \tE'\rightarrow \tE.
\end{equation}
On note 
\begin{equation}\label{higgs0-rmd10c}
\varphi_\ox\colon \tE\rightarrow X'_{\oeta,\fet}
\end{equation}
le foncteur composé $\beta'_*\circ\Phi^*$.

On désigne par $\fV_\ox$ la catégorie des $X$-schémas étales $\ox$-pointés, 
ou ce qui revient au même, la catégorie des voisinages de $\ox$ 
dans le site $\Et_{/X}$. Pour tout objet $(U,\xi\colon \ox\rightarrow U)$ de $\fV_\ox$, on note encore 
$\xi\colon X'\rightarrow U$ le $X$-morphisme induit par $\xi$. On démontre dans (\cite{ag2} 10.37)
que pour tout faisceau $F=\{U\mapsto F_U\}$ de $\tE$, on a un isomorphisme canonique fonctoriel 
\begin{equation}\label{higgs0-rmd10f}
\varphi_\ox(F)\stackrel{\sim}{\rightarrow}\underset{\underset{U\in \fV_\ox^\circ}{\longrightarrow}}{\lim}\ (\xi_\oeta)^*_\fet(F_U).
\end{equation}

Supposons que $\ox$ soit au-dessus de $s$. On démontre (\cite{ag3} 3.7) que $\oX'$ est normal et strictement local (et en particulier intègre). 
Soient $\oy$ un point géométrique de $X'_\oeta$ (qui est intègre), 
$\bB_{\pi_1(X'_\oeta,\oy)}$ le topos classifiant du groupe fondamental $\pi_1(X'_\oeta,\oy)$, 
\begin{equation}\label{higgs0-rmd10cd}
\nu_\oy\colon X'_{\oeta,\fet}\stackrel{\sim}{\rightarrow}\bB_{\pi_1(X'_\oeta,\oy)}
\end{equation}
le foncteur fibre en $\oy$ (\cite{ag2} (9.8.4)). 
Le foncteur composé
\begin{equation}\label{higgs0-rmd10d}
\xymatrix{
{\tE}\ar[r]^-(0.5){\varphi_\ox}&{X'_{\oeta,\fet}}\ar[r]^-(0.5){\nu_\oy}&{\bB_{\pi_1(X'_\oeta,\oy)}}\ar[r]&{\Ens}}
\end{equation} 
où la dernière flèche est le foncteur oubli de l'action de $\pi_1(X'_\oeta,\oy)$, est un foncteur fibre (\cite{ag2} 10.31 et 9.9). 
Il correspond à un point d'origine géométrique du topos $\tE$, noté $\rho(\oy\rightsquigarrow \ox)$  
(cf. \cite{ag3} 8.6).

\begin{teo}[\cite{ag2} 10.30]\label{higgs0-rmd100}
Sous les hypothèses de \eqref{higgs0-rmd10}, pour tout faisceau abélien $F$ de $\tE$ et tout entier $i\geq 0$, on a un isomorphisme canonique fonctoriel \eqref{higgs0-rmd7c}
\begin{equation}
\rR^i\sigma_*(F)_\ox\stackrel{\sim}{\rightarrow}\rH^i(X'_{\oeta,\fet},\varphi_\ox(F)).
\end{equation}
\end{teo}

\begin{cor}\label{higgs0-rmd101}
Sous les hypothèses de \eqref{higgs0-rmd10}, si de plus $\ox$ est au-dessus de $s$, 
pour tout faisceau abélien $F$ de $\tE$ et tout entier $i\geq 0$, on a un isomorphisme canonique fonctoriel 
\begin{equation}\label{higgs0-rmd10h}
\rR^i\sigma_*(F)_\ox\stackrel{\sim}{\rightarrow}\rH^i(\pi_1(X'_{\oeta},\oy),\nu_\oy(\varphi_\ox(F))).
\end{equation}
\end{cor}

\begin{prop}[\cite{ag2} 10.32] 
La famille des foncteurs $\varphi_\ox$ \eqref{higgs0-rmd10c},
lorsque $\ox$ décrit l'ensemble des points géométriques de $X$, est conservative.
\end{prop}

\subsection{}\label{higgs0-rmd2}
Pour chaque objet $(V\rightarrow U)$ de $E$, on note $\oU^V$ la clôture intégrale de $\oU=U\times_S\oS$
dans $V$ et on pose 
\begin{equation}\label{higgs0-rmd2a}
\ocB(V\rightarrow U)=\Gamma(\oU^V,\co_{\oU^V}).
\end{equation} 
On définit ainsi un préfaisceau d'anneaux sur $E$, qui se trouve être un {\em faisceau} (\cite{ag3} 8.16). 
On notera que $\ocB$ n'est pas en général un faisceau pour la topologie de $E$ définie originellement 
par Faltings dans (\cite{faltings2} page 214) (cf. \cite{ag3} 8.18). 
Pour tout $U\in \ob(\Et_{/X})$, on note $\ocB_U$ la restriction de $\ocB$ 
à la fibre $\Et_{\rf/U_\oeta}$ de $\pi$ au-dessus de $U$, de sorte que $\ocB=\{U\rightarrow \ocB_U\}$. 
On donne dans \ref{higgs0-rmd3} ci-dessous une description explicite de ce faisceau. Pour tout entier $n\geq 0$, on pose 
\begin{eqnarray}
\ocB_n&=&\ocB/p^n\ocB,\label{higgs0-rmd2c}\\
\ocB_{U,n}&=&\ocB_U/p^n\ocB_U.\label{higgs0-rmd2cc}
\end{eqnarray}
On notera que la correspondance $\{U\mapsto \ocB_{U,n}\}$
forme naturellement un préfaisceau sur $E$ dont le faisceau associé est canoniquement isomorphe à $\ocB_n$. 
Il est en général difficile, voir impossible, de décrire explicitement les restrictions de $\ocB_n$ 
aux fibres du foncteur $\pi$ \eqref{higgs0-rmd1b}. Toutefois, ses images par les foncteurs fibres 
\eqref{higgs0-rmd10d} sont accessibles (cf. \cite{ag3} (10.8.5)). 

On note $\hbar\colon \oX\rightarrow X$ la projection canonique \eqref{higgs0-not1a} et 
\begin{equation}\label{higgs0-rmd2d}
\hbar_*(\co_\oX)\rightarrow \sigma_*(\ocB)
\end{equation}
l'homomorphisme défini pour tout $U\in \ob(\Et_{/X})$ par l'homomorphisme canonique
\begin{equation}\label{higgs0-rmd2b}
\Gamma(\oU,\co_\oU)\rightarrow \Gamma(\oU^{U_\oeta},\co_{\oU^{U_\oeta}}).
\end{equation}
Sauf mention explicite du contraire, on considère $\sigma\colon \tE\rightarrow X_\et$ \eqref{higgs0-rmd7c}
comme un morphisme de topos annelés, respectivement par $\ocB$ et $\hbar_*(\co_\oX)$.

\subsection{}\label{higgs0-rmd3}
Soient $U$ un objet de $\Et_{/X}$, 
$\oy$ un point géométrique de $U_\oeta$, $V$ la composante connexe de $U_\oeta$ contenant $\oy$.
On désigne par $\bB_{\pi_1(V,\oy)}$ le topos classifiant du groupe fondamental $\pi_1(V,\oy)$, par $(V_i)_{i\in I}$ 
le revêtement universel normalisé de $V$ en $\oy$ (\cite{ag2} 9.8) et par 
\begin{equation}\label{higgs0-rmd3a}
\nu_\oy\colon V_\fet\stackrel{\sim}{\rightarrow} \bB_{\pi_1(V,\oy)},\ \ \ 
F\mapsto \underset{\underset{i\in I^\circ}{\longrightarrow}}{\lim}\ F(V_i)
\end{equation} 
le foncteur fibre en $\oy$. Pour tout $i\in I$, $(V_i\rightarrow U)$ est naturellement un objet de $E$.
On peut donc considérer le système projectif des schémas $(\oU^{V_i})_{i\in I}$. On pose  
\begin{equation}\label{higgs0-rmd3b}
\oR^\oy_U= \underset{\underset{i\in I^\circ}{\longrightarrow}}{\lim}\ \Gamma(\oU^{V_i},\co_{\oU^{V_i}}),
\end{equation} 
qui est un anneau de $\bB_{\pi_1(V,\oy)}$. D'après (\cite{ag3} 8.15), 
on a un isomorphisme canonique de $\bB_{\pi_1(V,\oy)}$ 
\begin{equation}\label{higgs0-rmd3c}
\nu_{\oy}(\ocB_U|V)\stackrel{\sim}{\rightarrow}\oR^\oy_U.
\end{equation}

\subsection{}\label{higgs0-rmd14}
Comme $X_\eta$ est un ouvert de $X_\et$ ({\em i.e.}, un sous-objet de l'objet final $X$), 
$\sigma^*(X_\eta)$ est un ouvert de $\tE$. On note 
\begin{equation}\label{higgs0-rmd14a}
\gamma\colon \tE_{/\sigma^*(X_\eta)}\rightarrow \tE
\end{equation}
le morphisme de localisation de $\tE$ en $\sigma^*(X_\eta)$.
On désigne par $\tE_s$ le sous-topos fermé de $\tE$ complémentaire de l'ouvert $\sigma^*(X_\eta)$, 
c'est-à-dire la sous-catégorie pleine de $\tE$ formée des faisceaux $F$ tels que $\gamma^*(F)$
soit un objet final de $\tE_{/\sigma^*(X_\eta)}$, et par 
\begin{equation}\label{higgs0-rmd14b}
\delta\colon \tE_s\rightarrow \tE
\end{equation} 
le plongement canonique, c'est-à-dire le morphisme de topos tel que  
$\delta_*\colon \tE_s\rightarrow \tE$ soit le foncteur d'injection canonique. Il existe un morphisme 
\begin{equation}\label{higgs0-TFT8c}
\sigma_s\colon \tE_s\rightarrow X_{s,\et}
\end{equation}
unique à isomorphisme près tel que le diagramme 
\begin{equation}\label{higgs0-TFT8d}
\xymatrix{
{\tE_s}\ar[r]^{\sigma_s}\ar[d]_{\delta}&{X_{s,\et}}\ar[d]^{\iota_\et}\\
{\tE}\ar[r]^\sigma&{X_\et}}
\end{equation}
où $\iota\colon X_s\rightarrow X$ est l'injection canonique,
soit commutatif à isomorphisme près (cf. \cite{ag3} 9.8).

Pour tout entier $n\geq 1$, identifiant les topos étales de $X_s$ et $\oX_n$  \eqref{higgs0-not1a}
($k$ étant algébriquement clos), le morphisme $\sigma_s$ et l'homomorphisme \eqref{higgs0-rmd2d} 
induisent un morphisme de topos annelés (\cite{ag3} 9.9)
\begin{equation}\label{higgs0-rmd14c}
\sigma_n\colon (\tE_s,\ocB_n)\rightarrow (X_{s,\et},\co_{\oX_n}).  
\end{equation}

\subsection{}\label{higgs0-rmd4}
Dans la suite de cette introduction, on suppose qu'il existe une $\cA_2(\oS)$-déformation lisse $\tX$ de $\coX$
que l'on fixe  (cf. \eqref{higgs0-not1a} et \ref{higgs0-dth1} pour les notations)~:
\begin{equation}
\xymatrix{
{\coX}\ar[r]\ar[d]&{\tX}\ar[d]\\
{\coS}\ar[r]&{\cA_2(\oS)}}
\end{equation}
Soit $Y=\Spec(R)$ un objet affine connexe de $\Et_{/X}$ admettant un $S$-morphisme étale 
vers $\mG_{m,S}^d$ pour un entier $d\geq 1$ et tel que $Y_s\not=\emptyset$
(autrement dit, $Y$ vérifie les hypothèses de \ref{higgs0-dth0}) et soit $\tY\rightarrow \tX$ l'unique morphisme étale 
qui relève $\coY\rightarrow \coX$. Pour tout point géométrique $\oy$ de $Y_\oeta$,  
on désigne par $Y'_\oeta$ la composante connexe de $Y_\oeta$ contenant $\oy$, par 
\begin{equation}
\nu_\oy\colon Y'_{\oeta,\fet}\stackrel{\sim}{\rightarrow}\bB_{\pi_1(Y'_\oeta,\oy)}
\end{equation} 
le foncteur fibre en $\oy$, 
et par $\cF^\oy_Y$ la $\hoR^\oy_Y$-extension de Higgs-Tate associée à $(Y,\tY)$ \eqref{higgs0-dth11}.

Pour tout entier $n\geq 0$, il existe un $\ocB_{Y,n}$-module $\cF_{Y,n}$, où $\ocB_{Y,n}=\ocB_Y/p^n\ocB_Y$ 
\eqref{higgs0-rmd2cc}, unique à isomorphisme canonique près, tel que 
pour tout point géométrique $\oy$ de $Y_\oeta$, on ait un isomorphisme canonique de $\oR^\oy_Y$-modules  \eqref{higgs0-rmd3c}
\begin{equation}\label{higgs0-rmd4e}
\nu_\oy(\cF_{Y,n}|Y'_\oeta)\stackrel{\sim}{\rightarrow}\cF^\oy_Y/p^n\cF^\oy_Y.
\end{equation}
La suite exacte \eqref{higgs0-dth11d} induit une suite exacte canonique de 
$\ocB_{Y,n}$-modules 
\begin{equation}\label{higgs0-RGG14a}
0\rightarrow \ocB_{Y,n}\rightarrow \cF_{Y,n}\rightarrow 
\xi^{-1}\Omega^1_{X/S}(Y)\otimes_{\co_X(Y)}\ocB_{Y,n} \rightarrow 0. 
\end{equation}

Pour tout nombre rationnel $r\geq 0$,    
on note $\cF^{(r)}_{Y,n}$ l'extension de $\ocB_{Y,n}$-modules de $Y_{\oeta,\fet}$ déduite de 
$\cF_{Y,n}$ par image inverse 
par le morphisme de multiplication par $p^r$ sur $\xi^{-1}\Omega^1_{X/S}(Y)\otimes_{\co_X(Y)}\ocB_{Y,n}$,
de sorte qu'on a une suite exacte canonique de $\ocB_{Y,n}$-modules
\begin{equation}\label{higgs0-RGG22a}
0\rightarrow \ocB_{Y,n}\rightarrow \cF^{(r)}_{Y,n}\rightarrow 
\xi^{-1}\Omega^1_{X/S}(Y)\otimes_{\co_X(Y)}\ocB_{Y,n} \rightarrow 0. 
\end{equation}
Celle-ci induit pour tout entier $m\geq 1$, une suite exacte de $\ocB_{Y,n}$-modules
\[
0\rightarrow \rS^{m-1}_{\ocB_{Y,n}}(\cF^{(r)}_{Y,n})\rightarrow \rS^{m}_{\ocB_{Y,n}}(\cF^{(r)}_{Y,n})\rightarrow 
\rS^m_{\ocB_{Y,n}}(\xi^{-1}\Omega^1_{X/S}(Y)\otimes_{\co_X(Y)} \ocB_{Y,n})\rightarrow 0.
\]
Les $\ocB_{Y,n}$-modules $(\rS^{m}_{\ocB_{Y,n}}(\cF^{(r)}_{Y,n}))_{m\in \mN}$ forment donc un système inductif 
dont la limite inductive 
\begin{equation}\label{higgs0-RGG22b}
\cC^{(r)}_{Y,n}=\underset{\underset{m\geq 0}{\longrightarrow}}\lim\ \rS^m_{\ocB_{Y,n}}(\cF^{(r)}_{Y,n})
\end{equation}
est naturellement munie d'une structure de $\ocB_{Y,n}$-algèbre de $Y_{\oeta,\fet}$. 

Pour tous nombres rationnels $r\geq r'\geq 0$, on a un morphisme $\ocB_{Y,n}$-linéaire canonique 
\begin{equation}\label{higgs0-RGG22g}
\tta_{Y,n}^{r,r'}\colon\cF_{Y,n}^{(r)}\rightarrow \cF_{Y,n}^{(r')}
\end{equation} 
qui relève la multiplication par $p^{r'-r}$ sur  
$\xi^{-1}\Omega^1_{X/S}(Y)\otimes_{\co_X(Y)}\ocB_{Y,n}$ et qui étend l'identité sur $\ocB_{Y,n}$ 
\eqref{higgs0-RGG22a}. Il induit un homomorphisme de $\ocB_{Y,n}$-algèbres 
\begin{equation}\label{higgs0-RGG22h}
\alpha_{Y,n}^{r,r'}\colon \cC_{Y,n}^{(r)}\rightarrow \cC_{Y,n}^{(r')}.
\end{equation}
On notera que $\cC^{(r)}_{Y,n}$ et $\cF^{(r)}_{Y,n}$ dépendent du choix de la déformation $\tX$.

On étend les définitions précédentes aux objets affines connexes $Y$ de $\Et_{/X}$ 
tels que $Y_s=\emptyset$ en posant $\cC^{(r)}_{Y,n}=\cF^{(r)}_{Y,n}=0$.

\subsection{}\label{higgs0-rmd6}
Soient $n$ entier $\geq 0$, $r$ un nombre rationnel $\geq 0$. 
Les correspondances $\{Y\mapsto \cF^{(r)}_{Y,n}\}$ et $\{Y\mapsto \cC^{(r)}_{Y,n}\}$ 
forment naturellement des préfaisceaux sur la sous-catégorie pleine de $E$ formée des objets $(V\rightarrow Y)$
tels que $Y$ soit affine, connexe et admette un morphisme étale 
vers $\mG_{m,S}^d$ pour un entier $d\geq 1$ (\cite{ag3} 10.19). Cette sous-catégorie étant clairement topologiquement
génératrice dans $E$, on peut considérer les faisceaux associés dans $\tE$
\begin{eqnarray}
\cF^{(r)}_n&=&\{Y\mapsto \cF^{(r)}_{Y,n}\}^a,\label{higgs0-rmd6b}\\
\cC^{(r)}_n&=&\{Y\mapsto \cC^{(r)}_{Y,n}\}^a. \label{higgs0-rmd6a}
\end{eqnarray}
Ce sont en fait un $\ocB_n$-module et une $\ocB_n$-algèbre de $\tE_s$ (\cite{ag3} 10.22). 
On appelle $\cF^{(r)}_n$ la {\em $\ocB_n$-extension de Higgs-Tate d'épaisseur $r$} 
et $\cC^{(r)}_n$ {\em la $\ocB_n$-algèbre de Higgs-Tate d'épaisseur $r$} associés à $\tX$. 
On a une suite exacte canonique de $\ocB_n$-modules \eqref{higgs0-rmd14c}
\begin{equation}\label{higgs0-rmd6c}
0\rightarrow \ocB_n\rightarrow \cF^{(r)}_n\rightarrow 
\sigma_n^*(\xi^{-1}\Omega^1_{\oX_n/\oS_n})\rightarrow 0.
\end{equation}
On décrit explicitement dans (\cite{ag3} 10.29) les images de $\cF_n^{(r)}$ et $\cC_n^{(r)}$ 
par les foncteurs fibres \eqref{higgs0-rmd10c}. 

Pour tous nombres rationnels $r\geq r'\geq 0$, les homomorphismes \eqref{higgs0-RGG22h} induisent un homomorphisme de 
$\ocB_n$-algèbres 
\begin{equation}\label{higgs0-RGG24d}
\alpha_n^{r,r'}\colon \cC_n^{(r)}\rightarrow \cC_n^{(r')}.
\end{equation}
Pour tous nombres rationnels $r\geq r'\geq r''\geq 0$, on a
\begin{equation}\label{higgs0-RGG18h}
\alpha_n^{r,r''}=\alpha_n^{r',r''} \circ \alpha_n^{r,r'}.
\end{equation}

On a un isomorphisme $\cC_n^{(r)}$-linéaire canonique 
\begin{equation}\label{higgs0-rmd6i}
\Omega^1_{\cC_n^{(r)}/\ocB_n}\stackrel{\sim}{\rightarrow} 
\sigma_n^*(\xi^{-1}\Omega^1_{\oX_n/\oS_n})\otimes_{\ocB_n}\cC_n^{(r)}.
\end{equation}
La $\ocB_n$-dérivation universelle de $\cC_n^{(r)}$ correspond via cet isomorphisme à l'unique $\ocB_n$-dérivation 
\begin{equation}\label{higgs0-rmd6j}
d_n^{(r)}\colon \cC_n^{(r)}\rightarrow \sigma_n^*(\xi^{-1}\Omega^1_{\oX_n/\oS_n})\otimes_{\ocB_n}\cC_n^{(r)}
\end{equation}
qui prolonge le morphisme canonique $\cF_n^{(r)}\rightarrow 
\sigma_n^*(\xi^{-1}\Omega^1_{\oX_n/\oS_n})$.
Pour tous nombres rationnels $r\geq r'\geq 0$, on a 
\begin{equation}\label{higgs0-rmd6l}
p^{r-r'}(\id \otimes \alpha^{r,r'}_n) \circ d^{(r)}_n=d^{(r')}_n\circ \alpha^{r,r'}_n.
\end{equation}

\subsection{}\label{higgs0-rmd9}
Les systèmes projectifs d'objets de $\tE_s$ \eqref{higgs0-rmd7} indexés par l'ensemble ordonné 
des entiers naturels $\mN$, forment un topos que l'on note $\tE_s^{\mN^\circ}$. On trouvera dans 
(\cite{ag3} § 7) des sorites utiles sur ce type de topos. On désigne par $\bvocB$ l'anneau 
$(\ocB_{n+1})_{n\in \mN}$ de $\tE_s^{\mN^\circ}$, 
par $\co_{\bvoX}$ l'anneau $(\co_{\oX_{n+1}})_{n\in \mN}$ de $X_{s,\et}^{\mN^\circ}$,
et par $\xi^{-1}\Omega^1_{\bvoX/\bvoS}$ le $\co_{\bvoX}$-module 
$(\xi^{-1}\Omega^1_{\oX_{n+1}/\oS_{n+1}})_{n\in \mN}$ de $X_{s,\et}^{\mN^\circ}$. 
Les morphismes $(\sigma_{n+1})_{n\in \mN}$ \eqref{higgs0-rmd14c} induisent un morphisme de topos annelés
\begin{equation}\label{higgs0-rmd9a}
\bvsigma\colon (\tE_s^{\mN^\circ},\bvocB)\rightarrow(X_{s,\et}^{\mN^\circ},\co_{\bvoX}). 
\end{equation}

On dit qu'un $\bvocB$-module $(M_{n+1})_{n\in \mN}$
de $\tE_s^{\mN^\circ}$ est {\em adique} si pour tous entiers $m$ et $n$ tels que $m\geq n\geq 1$, 
le morphisme $M_m\otimes_{\ocB_m}\ocB_n\rightarrow M_n$ déduit du morphisme de transition 
$M_m\rightarrow M_n$ est un isomorphisme.

Soit $r$ un nombre rationnel $\geq 0$. Pour tous entiers $m\geq n\geq 1$, 
on a un morphisme $\ocB_m$-linéaire canonique
$\cF_m^{(r)}\rightarrow \cF_n^{(r)}$ et un homomorphisme canonique de $\ocB_m$-algèbres 
$\cC_m^{(r)}\rightarrow \cC_n^{(r)}$ tels que les morphismes induits 
\begin{equation}\label{higgs0-rmd9b}
\cF_m^{(r)}\otimes_{\ocB_m}\ocB_n\rightarrow \cF_n^{(r)}\ \ \ 
{\rm et}\ \ \ \cC_m^{(r)}\otimes_{\ocB_m}\ocB_n\rightarrow \cC_n^{(r)}
\end{equation}
soient des isomorphismes. Ces morphismes forment des systèmes compatibles lorsque $m$ et $n$ varient, 
de sorte que $(\cF_{n+1}^{(r)})_{n\in \mN}$ et $(\cC_{n+1}^{(r)})_{n\in \mN}$ sont des systèmes projectifs. 
On appelle {\em $\bvocB$-extension de Higgs-Tate d'épaisseur $r$} associée à $\tX$, 
et l'on note $\bvcF^{(r)}$, le $\bvocB$-module $(\cF_{n+1}^{(r)})_{n\in \mN}$ de $\tE_s^{\mN^\circ}$. 
On appelle {\em $\bvocB$-algèbre de Higgs-Tate d'épaisseur $r$} associée à $\tX$,
et l'on note $\bvcC^{(r)}$, la $\bvocB$-algèbre $(\cC_{n+1}^{(r)})_{n\in \mN}$ de $\tE_s^{\mN^\circ}$.
Ce sont des $\bvocB$-modules adiques. 
On a une suite exacte de $\bvocB$-modules 
\begin{equation}\label{higgs0-rmd9c}
0\rightarrow \bvocB\rightarrow \bvcF^{(r)}\rightarrow 
\bvsigma^*(\xi^{-1}\Omega^1_{\bvoX/\bvoS})\rightarrow 0.
\end{equation}
Pour tous nombres rationnels $r\geq r'\geq 0$, les homomorphismes $(\alpha_n^{r,r'})_{n\in \mN}$
induisent un homomorphisme de $\bvocB$-algèbres 
\begin{equation}\label{higgs0-rmd9e}
\bvalpha^{r,r'}\colon \bvcC^{(r)}\rightarrow \bvcC^{(r')}.
\end{equation}
Pour tous nombres rationnels $r\geq r'\geq r''\geq 0$, on a
\begin{equation}\label{higgs0-rmd9g}
\bvalpha^{r,r''}=\bvalpha^{r',r''} \circ \bvalpha^{r,r'}.
\end{equation}

Les dérivations $(d_{n+1}^{(r)})_{n\in \mN}$ définissent un morphisme
\begin{equation}\label{higgs0-rmd9d}
\bvd^{(r)}\colon \bvcC^{(r)}\rightarrow \bvsigma^*(\xi^{-1}\Omega^1_{\bvoX/\bvoS})\otimes_{\bvocB}\bvcC^{(r)},
\end{equation}
qui n'est autre que la $\bvocB$-dérivation universelle de $\bvcC^{(r)}$. 
Pour tous nombres rationnels $r\geq r'\geq 0$, on a 
\begin{equation}\label{higgs0-rmd9f}
p^{r-r'}(\id \otimes \bvalpha^{r,r'}) \circ \bvd^{(r)}=\bvd^{(r')}\circ \bvalpha^{r,r'}.
\end{equation}

\section{Modules de Dolbeault}\label{higgs0-md}

\subsection{}\label{higgs0-md1}
Les hypothèses et notations de § \ref{higgs0-rmd} sont en vigueur dans cette section. 
On pose $\cS=\Spf(\co_C)$ et on note $\fX$ le schéma formel complété $p$-adique de $\oX$,
et $\xi^{-1}\Omega^1_{\fX/\cS}$ le complété $p$-adique du $\co_\oX$-module 
$\xi^{-1}\Omega^1_{\oX/\oS}=\xi^{-1}\Omega^1_{X/S}\otimes_{\co_X}\co_{\oX}$.
On désigne par
\begin{equation}\label{higgs0-md1a}
\bvu\colon (X_{s,\et}^{\mN^\circ},\co_{\bvoX})\rightarrow (X_{s,\zar}^{\mN^\circ},\co_{\bvoX})
\end{equation}
le morphisme canonique de topos annelés (\ref{higgs0-rmd9} et \cite{ag3} 2.9), par 
\begin{equation}\label{higgs0-md1b}
\lambda\colon (X_{s,\zar}^{\mN^\circ},\co_{\bvoX})\rightarrow (X_{s,\zar}, \co_\fX)
\end{equation}
le morphisme de topos annelés dont le foncteur image directe correspondant est la limite projective (\cite{ag3} 7.4), 
et par
\begin{equation}\label{higgs0-md1c}
\top\colon (\tE_s^{\mN^\circ},\bvocB)\rightarrow (\fX_{\zar},\co_{\fX})
\end{equation}
le morphisme composé   $\lambda\circ \bvu\circ \bvsigma$ \eqref{higgs0-rmd9a}. 
Nous utilisons pour les modules la notation $\top^{-1}$ pour désigner l'image
inverse au sens des faisceaux abéliens et nous réservons la notation 
$\top^*$ pour l'image inverse au sens des modules~; de même pour $\bvsigma$. 

On note
\begin{equation}\label{higgs0-md1d}
\delta\colon \xi^{-1}\Omega^1_{\fX/\cS}\rightarrow \rR^1\top_*(\bvocB)
\end{equation}
le morphisme $\co_{\fX}$-linéaire de $X_{s,\zar}$ composé du morphisme d'adjonction (cf. \cite{ag3} (11.2.5))
\begin{equation}\label{higgs0-md1e}
\xi^{-1}\Omega^1_{\fX/\cS}\rightarrow \top_*(\bvsigma^*(\xi^{-1}\Omega^1_{\bvoX/\bvoS}))
\end{equation}
et du morphisme bord de la suite exacte longue de cohomologie déduite
de la suite exacte canonique 
\begin{equation}\label{higgs0-md1f}
0\rightarrow \bvocB\rightarrow \bvcF\rightarrow 
\bvsigma^*(\xi^{-1}\Omega^1_{\bvoX/\bvoS})\rightarrow 0.
\end{equation}
On notera que le morphisme 
\begin{equation}\label{higgs0-md1g}
\top^*(\xi^{-1}\Omega^1_{\fX/\cS})\rightarrow \bvsigma^*(\xi^{-1}\Omega^1_{\bvoX/\bvoS})
\end{equation}
adjoint  de \eqref{higgs0-md1e} est un isomorphisme (\cite{ag3} (11.2.6)). 

\begin{teo}[\cite{ag3} 11.8]\label{higgs0-md2}
Il existe un et un unique isomorphisme de $\co_\fX[\frac 1 p]$-algèbres graduées 
\begin{equation}\label{higgs0-md2a}
\wedge(\xi^{-1}\Omega^1_{\fX/\cS}[\frac 1 p])\stackrel{\sim}{\rightarrow} \oplus_{i\geq 0}\rR^i\top_*(\bvocB)[\frac 1 p]
\end{equation}
dont la composante en degré un est le morphisme $\delta\otimes_{\mZ_p}\mQ_p$ \eqref{higgs0-md1d}.
\end{teo} 

Cet énoncé est la clé de voûte de l'approche de Faltings en théorie de Hodge $p$-adique. 
On le retrouve ici et là sous différentes incarnations.  Sa forme galoisienne locale (\cite{ag1} 8.21) 
est une conséquence du théorème de presque pureté de Faltings (\cite{ag1} 6.16). 
L'énoncé global a une variante entière (\cite{ag3} 11.3) qui se déduit du cas local 
par localisation \eqref{higgs0-rmd101}. 

\subsection{}\label{higgs0-md3}
La suite exacte canonique \eqref{higgs0-md1f} induit, pour tout entier $m\geq 1$, une suite exacte 
\begin{equation}\label{higgs0-md3a}
0\rightarrow \Sym^{m-1}_{\bvocB}(\bvcF)\rightarrow \Sym^m_{\bvocB}(\bvcF)\rightarrow 
\bvsigma^*(\Sym^m_{\co_{\bvoX}}(\xi^{-1}\Omega^1_{\bvoX/\bvoS}))\rightarrow 0.
\end{equation}

\begin{prop}[\cite{ag3} 11.12]\label{higgs0-md4}
Soit $m$ un entier $\geq 1$.  
Alors~:
\begin{itemize}
\item[{\rm (i)}] Le morphisme  
\begin{equation}\label{higgs0-md4a}
\top_*(\Sym^{m-1}_{\bvocB}(\bvcF))[\frac 1 p]\rightarrow \top_*(\Sym^m_{\bvocB}(\bvcF))[\frac 1 p]
\end{equation}
induit par \eqref{higgs0-md3a} est un isomorphisme.
\item[{\rm (ii)}] Pour tout entier $q\geq 1$, le morphisme 
\begin{equation}\label{higgs0-md4b}
\rR^q\top_*(\Sym^{m-1}_{\bvocB}(\bvcF))[\frac 1 p]\rightarrow \rR^q\top_*(\Sym^m_{\bvocB}(\bvcF))[\frac 1 p]
\end{equation}
induit par \eqref{higgs0-md3a} est nul.
\end{itemize}
\end{prop}

La variante galoisienne locale de cet énoncé est due à Hyodo (\cite{hyodo1} 1.2). 
C'est le principal ingrédient dans la définition des systèmes locaux de Hodge-Tate. 

\begin{prop}[\cite{ag3} 11.18]\label{higgs0-md5}
L'homomorphisme canonique 
\begin{equation}
\co_{\fX}[\frac 1 p]\rightarrow \underset{\underset{r\in \mQ_{>0}}{\longrightarrow}}{\lim}\  \top_*(\bvcC^{(r)})[\frac 1 p]
\end{equation}
est un isomorphisme, et pour tout entier $q\geq 1$,
\begin{equation}
\underset{\underset{r\in \mQ_{>0}}{\longrightarrow}}{\lim}\ \rR^q\top_*(\bvcC^{(r)})[\frac 1 p] =0.
\end{equation}
\end{prop}

La variante galoisienne locale de cet énoncé (\cite{ag1} 12.5) est essentiellement due à Tsuji (\cite{tsuji3} 5.3.4). 

\subsection{}\label{higgs0-md8}
On note $\bMod(\bvocB)$ la catégorie des $\bvocB$-modules de $\tE_s^{\mN^\circ}$,
$\bMod^\ad(\bvocB)$ (resp. $\bMod^\atf(\bvocB)$) la sous-catégorie pleine formée des $\bvocB$-modules adiques 
(resp. adiques de type fini) \eqref{higgs0-rmd9} et $\bMod_{\mQ}(\bvocB)$ 
(resp. $\bMod^\ad_{\mQ}(\bvocB)$, resp. $\bMod^\atf_{\mQ}(\bvocB)$) 
la catégorie des objets de $\bMod(\bvocB)$ (resp. $\bMod^\ad(\bvocB)$, resp. $\bMod^\atf(\bvocB)$) 
à isogénie près \eqref{higgs0-not3}. La catégorie $\bMod_\mQ(\bvocB)$ est abélienne et les foncteurs canoniques 
\begin{equation}\label{higgs0-md8a}
\bMod^\atf_{\mQ}(\bvocB)\rightarrow \bMod^\ad_{\mQ}(\bvocB)\rightarrow\bMod_{\mQ}(\bvocB)
\end{equation} 
sont pleinement fidèles. 
On désigne par $\bMod^\coh(\co_{\fX})$ (resp. $\bMod^{\coh}(\co_{\fX}[\frac 1 p])$)
la catégorie des $\co_{\fX}$-modules (resp. $\co_{\fX}[\frac 1 p]$-modules) cohérents de $X_{s,\zar}$
et par $\bMod^{\coh}_\mQ(\co_{\fX})$ la catégorie des $\co_\fX$-modules cohérents à isogénie près. 
D'après (\cite{ag3} 6.16), le foncteur canonique
\begin{equation}\label{higgs0-md8b}
\bMod^{\coh}(\co_{\fX})\rightarrow \bMod^{\coh}(\co_{\fX}[\frac 1 p]),\ \ \ \cF\mapsto \cF_{\mQ_p}=\cF\otimes_{\mZ_p}\mQ_p,
\end{equation}
induit une équivalence de catégories abéliennes
\begin{equation}\label{higgs0-md8c}
\bMod^{\coh}_\mQ(\co_{\fX})\stackrel{\sim}{\rightarrow} \bMod^{\coh}(\co_{\fX}[\frac 1 p]). 
\end{equation}

\subsection{}\label{higgs0-md6}
Pour tout nombre rationnel $r\geq 0$, on note encore 
\begin{equation}\label{higgs0-md6a}
\bvd^{(r)}\colon \bvcC^{(r)}\rightarrow \top^*(\xi^{-1}\Omega^1_{\fX/\cS})\otimes_{\bvocB}\bvcC^{(r)}
\end{equation}
la $\bvocB$-dérivation induite par $\bvd^{(r)}$ \eqref{higgs0-rmd9d} et l'isomorphisme \eqref{higgs0-md1g},
que l'on identifie à la $\bvocB$-dérivation universelle de $\bvcC^{(r)}$.
C'est un $\bvocB$-champ de Higgs à coefficients dans $\top^*(\xi^{-1}\Omega^1_{\fX/\cS})$ \eqref{higgs0-not7}. 
On désigne par $\mK^\bullet(\bvcC^{(r)},p^r\bvd^{(r)})$ le complexe de Dolbeault
de $(\bvcC^{(r)},p^r\bvd^{(r)})$ \eqref{higgs0-not2} et par $\mK^\bullet_\mQ(\bvcC^{(r)},p^r\bvd^{(r)})$ 
son image dans $\bMod_{\mQ}(\bvocB)$. D'après \eqref{higgs0-rmd9f}, pour tous nombres rationnels $r\geq r'\geq 0$, 
l'homomorphisme $\bvalpha^{r,r'}$ \eqref{higgs0-rmd9e} induit un morphisme de complexes 
\begin{equation}\label{higgs0-md6b}
\bvnu^{r,r'}\colon \mK^\bullet(\bvcC^{(r)},p^r\bvd^{(r)})\rightarrow \mK^\bullet(\bvcC^{(r')},p^{r'}\bvd^{(r')}).
\end{equation}

\begin{prop}[\cite{ag3} 11.24]\label{higgs0-md7}
Le morphisme canonique 
\begin{equation}\label{higgs0-md7a}
\bvocB_\mQ\rightarrow \underset{\underset{r\in \mQ_{>0}}{\longrightarrow}}{\lim}\ 
\rH^0(\mK^\bullet_\mQ(\bvcC^{(r)},p^r\bvd^{(r)}))
\end{equation}
est un isomorphisme, et pour tout entier $q\geq 1$, 
\begin{equation}\label{higgs0-md7b}
\underset{\underset{r\in \mQ_{>0}}{\longrightarrow}}{\lim}\ \rH^q(\mK^\bullet_\mQ(\bvcC^{(r)},p^r\bvd^{(r)}))=0.
\end{equation}
\end{prop}

On notera que les limites inductives filtrantes ne sont pas a priori représentables dans la catégorie $\bMod_\mQ(\bvocB)$. 

\subsection{}\label{higgs0-md9}
Le foncteur $\top_*$ \eqref{higgs0-md1c} induit un foncteur additif et exact à gauche que l'on note encore 
\begin{equation}\label{higgs0-md9a}
\top_*\colon \bMod_{\mQ}(\bvocB) \rightarrow \bMod(\co_{\fX}[\frac 1 p]).
\end{equation}
Pour tout entier $q\geq 0$, on note 
\begin{equation}\label{higgs0-md9ab}
\rR^q \top_*\colon \bMod_\mQ(\bvocB)\rightarrow \bMod(\co_{\fX}[\frac 1 p])
\end{equation}
le $q$-ième foncteur dérivé droit de $\top_*$. 
D'après \eqref{higgs0-md8c}, le foncteur image inverse $\top^*$ induit un foncteur additif que l'on note encore
\begin{equation}\label{higgs0-md9b}
\top^*\colon \bMod^\coh(\co_{\fX}[\frac 1 p]) \rightarrow \bMod_{\mQ}^\atf(\bvocB). 
\end{equation}
Pour tout $\co_\fX[\frac 1 p]$-module cohérent $\cF$ et tout $\bvocB_\mQ$-module $\cG$, 
on a un homomorphisme canonique bifonctoriel 
\begin{equation}\label{higgs0-md9c}
\Hom_{\bvocB_\mQ}(\top^*(\cF),\cG)\rightarrow\Hom_{\co_\fX[\frac 1 p]}(\cF,\top_*(\cG)).
\end{equation}
qui est injectif (\cite{ag3} (12.1.5)).

\subsection{}\label{higgs0-md10}
On désigne par $\bIH(\co_\fX,\xi^{-1}\Omega^1_{\fX/\cS})$ la catégorie des $\co_\fX$-isogénies  
de Higgs à coefficients dans $\xi^{-1}\Omega^1_{\fX/\cS}$ \eqref{higgs0-imh1} et
par $\bIH^\coh(\co_\fX,\xi^{-1}\Omega^1_{\fX/\cS})$ la sous-catégorie pleine 
formée des quadruplets $(\cM,\cN,u,\theta)$ tels que $\cM$ et $\cN$ soient des $\co_\fX$-modules cohérents.
Ce sont des catégories additives. 
On note $\bIH_\mQ(\co_\fX,\xi^{-1}\Omega^1_{\fX/\cS})$ (resp. $\bIH^\coh_\mQ(\co_\fX,\xi^{-1}\Omega^1_{\fX/\cS})$)
la catégorie des objets de $\bIH(\co_\fX,\xi^{-1}\Omega^1_{\fX/\cS})$ 
(resp. $\bIH^\coh(\co_\fX,\xi^{-1}\Omega^1_{\fX/\cS})$) à isogénie près  \eqref{higgs0-not3}. 

On sous-entend par $\co_\fX[\frac 1 p]$-module de Higgs à coefficients dans $\xi^{-1}\Omega^1_{\fX/\cS}$, 
un $\co_\fX[\frac 1 p]$-module de Higgs à coefficients dans $\xi^{-1}\Omega^1_{\fX/\cS}[\frac 1 p]$ \eqref{higgs0-not2}. 
On désigne par $\bMH(\co_\fX[\frac 1 p], \xi^{-1}\Omega^1_{\fX/\cS})$ la catégorie de tels modules  
et par $\bMH^\coh(\co_\fX[\frac 1 p], \xi^{-1}\Omega^1_{\fX/\cS})$
la sous-catégorie pleine formée des modules
de Higgs dont le $\co_\fX[\frac 1 p]$-module sous-jacent est cohérent.

Le foncteur 
\begin{equation}\label{higgs0-MF12a}
\begin{array}[t]{clcr}
\bIH(\co_\fX,\xi^{-1}\Omega^1_{\fX/\cS})&\rightarrow& \bMH(\co_\fX[\frac 1 p], \xi^{-1}\Omega^1_{\fX/\cS})\\
(\cM,\cN,u,\theta)&\mapsto& (\cM_{\mQ_p}, (\id \otimes u_{\mQ_p}^{-1})\circ\theta_{\mQ_p})
\end{array}
\end{equation}
induit un foncteur 
\begin{equation}\label{higgs0-MF12b}
\bIH_\mQ(\co_\fX,\xi^{-1}\Omega^1_{\fX/\cS})\rightarrow \bMH(\co_\fX[\frac 1 p], \xi^{-1}\Omega^1_{\fX/\cS}).
\end{equation}
D'après (\cite{ag3} 6.21), celui-ci induit une équivalence de catégories 
\begin{equation}\label{higgs0-MF12c}
\bIH^\coh_\mQ(\co_\fX,\xi^{-1}\Omega^1_{\fX/\cS})\stackrel{\sim}{\rightarrow} 
\bMH^\coh(\co_\fX[\frac 1 p], \xi^{-1}\Omega^1_{\fX/\cS}).
\end{equation}

\begin{defi}\label{higgs0-md17}
On appelle {\em $\co_\fX[\frac 1 p]$-fibré de Higgs à coefficients dans $\xi^{-1}\Omega^1_{\fX/\cS}$} 
tout $\co_\fX[\frac 1 p]$-module de Higgs à coefficients dans $\xi^{-1}\Omega^1_{\fX/\cS}$ 
dont le $\co_\fX[\frac 1 p]$-module sous-jacent est localement projectif de type fini (\cite{ag3} 2.8). 
\end{defi}

\subsection{}\label{higgs0-md11}
Soit $r$ un nombre rationnel $\geq 0$. On désigne par $\Xi^r$ la catégorie des $p^r$-isoconnexions intégrables 
relativement à l'extension $\bvcC^{(r)}/\bvocB$ \eqref{higgs0-not5} et par $\Xi^r_\mQ$ la catégorie des objets de 
$\Xi^r$ à isogénie près  \eqref{higgs0-not3}. On note $\fS^r$ le foncteur 
\begin{equation}\label{higgs0-md11a}
\fS^r\colon \bMod(\bvocB)\rightarrow \Xi^r, \ \ \ \cM\mapsto (\bvcC^{(r)}\otimes_{\bvocB}\cM,\bvcC^{(r)}\otimes_{\bvocB}\cM,
\id,p^r\bvd^{(r)}\otimes \id).
\end{equation}
Celui-ci induit un foncteur que l'on note encore 
\begin{equation}\label{higgs0-md11aa}
\fS^r\colon \bMod_\mQ(\bvocB)\rightarrow \Xi^r_\mQ.
\end{equation}
On note $\cK^r$ le foncteur 
\begin{equation}\label{higgs0-md11b}
\cK^r\colon \Xi^r\rightarrow \bMod(\bvocB),\ \ \ (\cF,\cG,u,\nabla)\mapsto \ker(\nabla).
\end{equation}
Celui-ci induit un foncteur que l'on note encore 
\begin{equation}\label{higgs0-md11bb}
\cK^r\colon \Xi^r_\mQ\rightarrow \bMod_\mQ(\bvocB).
\end{equation}
Il est clair que \eqref{higgs0-md11a} est un adjoint à gauche de \eqref{higgs0-md11b}. 
Par suite, \eqref{higgs0-md11aa} est un adjoint à gauche de \eqref{higgs0-md11bb}.

Si $(\cN,\cN',v,\theta)$ est une $\co_\fX$-isogénie de Higgs
à coefficients dans $\xi^{-1}\Omega^1_{\fX/\cS}$,
\begin{equation}\label{higgs0-md11c}
(\bvcC^{(r)}\otimes_{\bvocB}\top^*(\cN),\bvcC^{(r)}\otimes_{\bvocB}\top^*(\cN'),\id \otimes_{\bvocB}\top^*(v),
\id \otimes \top^*(\theta)+p^r \bvd^{(r)} \otimes\top^*(v))
\end{equation}
est un objet de $\Xi^r$ (\cite{ag3} 6.12). On obtient ainsi un foncteur \eqref{higgs0-md10}
\begin{equation}\label{higgs0-md11d}
\top^{r+}\colon \bIH(\co_\fX,\xi^{-1}\Omega^1_{\fX/\cS}) \rightarrow \Xi^r.
\end{equation}
D'après \eqref{higgs0-MF12c}, celui-ci induit un foncteur que l'on note encore
\begin{equation}\label{higgs0-md11e}
\top^{r+}\colon\bMH^\coh(\co_\fX[\frac 1 p], \xi^{-1}\Omega^1_{\fX/\cS})\rightarrow \Xi^r_\mQ.
\end{equation}

Soit $(\cF,\cG,u,\nabla)$ un objet de $\Xi^r$. D'après la formule de projection (\cite{ag3} 12.4), $\nabla$ 
induit un morphisme $\co_\fX$-linéaire~: 
\begin{equation}\label{higgs0-md11f}
\top_*(\nabla)\colon \top_*(\cF)\rightarrow \xi^{-1}\Omega^1_{\fX/\cS}\otimes_{\co_\fX}\top_*(\cG).
\end{equation}
On voit aussitôt que $(\top_*(\cF),\top_*(\cG),\top_*(u),\top_*(\nabla))$ est une $\co_\fX$-isogénie 
de Higgs à coefficients dans $\xi^{-1}\Omega^1_{\fX/\cS}$. On obtient ainsi un foncteur
\begin{equation}\label{higgs0-md11g}
\top^r_+\colon \Xi^r\rightarrow \bIH(\co_\fX,\xi^{-1}\Omega^1_{\fX/\cS}),
\end{equation}
qui est clairement un adjoint à droite de \eqref{higgs0-md11d}. 
Le composé des foncteurs \eqref{higgs0-md11g} et \eqref{higgs0-MF12a} induit un foncteur  que l'on note encore
\begin{equation}\label{higgs0-md11h}
\top^r_+\colon \Xi^r_\mQ\rightarrow \bMH(\co_\fX[\frac 1 p], \xi^{-1}\Omega^1_{\fX/\cS}).
\end{equation}

\begin{defi}[\cite{ag3} 12.10]\label{higgs0-md13}
Soient $\cM$ un objet de $\bMod^\atf_{\mQ}(\bvocB)$, 
$\cN$ un $\co_\fX[\frac 1 p]$-fibré de Higgs à coefficients dans $\xi^{-1}\Omega^1_{\fX/\cS}$.
\begin{itemize}
\item[(i)] Soit $r$ un nombre rationnel $>0$.  On dit que $\cM$ et  $\cN$ sont {\em $r$-associés} 
s'il existe un isomorphisme de $\Xi^r_\mQ$ 
\begin{equation}\label{higgs0-md13a}
\alpha\colon \top^{r+}(\cN) \stackrel{\sim}{\rightarrow}\fS^r(\cM).
\end{equation}
On dit alors aussi que le triplet $(\cM,\cN,\alpha)$ est {\em $r$-admissible}. 
\item[(ii)] On dit que $\cM$ et  $\cN$ sont {\em associés} s'il existe un nombre rationnel $r>0$ tel que 
$\cM$ et  $\cN$ soient $r$-associés.
\end{itemize}
\end{defi}

On notera que pour tous nombres rationnels $r\geq r'>0$, 
si $\cM$ et  $\cN$ sont $r$-associés, ils sont $r'$-associés. 

\begin{defi}[\cite{ag3} 12.11] \label{higgs0-md14}
(i)\ On appelle {\em $\bvocB_\mQ$-module de Dolbeault} tout objet de $\bMod^\atf_{\mQ}(\bvocB)$
pour lequel il existe un $\co_\fX[\frac 1 p]$-fibré de Higgs associé à coefficients dans $\xi^{-1}\Omega^1_{\fX/\cS}$.

(ii)\ On dit qu'un $\co_\fX[\frac 1 p]$-fibré de Higgs à coefficients dans $\xi^{-1}\Omega^1_{\fX/\cS}$ est {\em soluble}
s'il admet un module de Dolbeault associé. 
\end{defi}

On désigne par $\bMod_\mQ^\Dolb(\bvocB)$ la sous-catégorie pleine de $\bMod^\atf_\mQ(\bvocB)$ 
formée des $\bvocB_\mQ$-modules de Dolbeault, et par $\bMH^\sol(\co_\fX[\frac 1 p], \xi^{-1}\Omega^1_{\fX/\cS})$ 
la sous-catégorie pleine de  $\bMH(\co_\fX[\frac 1 p], \xi^{-1}\Omega^1_{\fX/\cS})$
formée des $\co_\fX[\frac 1 p]$-fibrés de Higgs solubles à coefficients dans $\xi^{-1}\Omega^1_{\fX/\cS}$.

\subsection{}\label{higgs0-md12}
Pour tout $\bvocB_\mQ$-module $\cM$ et tous nombres rationnels $r\geq r'\geq 0$, 
on a un morphisme canonique de $\bMH(\co_\fX[\frac 1 p], \xi^{-1}\Omega^1_{\fX/\cS})$ 
\begin{equation}\label{higgs0-md12a}
\top^r_+(\fS^r(\cM))\rightarrow \top^{r'}_+(\fS^{r'}(\cM)).
\end{equation}
On obtient ainsi un système inductif filtrant $(\top^r_+(\fS^r(\cM)))_{r\in \mQ_{\geq 0}}$. 
On désigne par $\cH$ le foncteur 
\begin{equation}\label{higgs0-md12b}
\cH\colon \bMod_\mQ(\bvocB)\rightarrow \bMH(\co_\fX[\frac 1 p], \xi^{-1}\Omega^1_{\fX/\cS}), \ \ \ \cM\mapsto 
\underset{\underset{r\in \mQ_{>0}}{\longrightarrow}}{\lim}\ \top^r_+(\fS^r(\cM)). 
\end{equation}

Pour tout objet $\cN$ de $\bMH(\co_\fX[\frac 1 p], \xi^{-1}\Omega^1_{\fX/\cS})$ 
et tous nombres rationnels $r\geq r'\geq 0$, on a un morphisme canonique  de $\bMod_{\mQ}(\bvocB)$
\begin{equation}\label{higgs0-md12c}
\cK^r(\top^{r+}(\cN))\rightarrow \cK^{r'}(\top^{r'+}(\cN)).
\end{equation}
On obtient ainsi un système inductif filtrant $(\cK^r(\top^{r+}(\cN)))_{r\geq  0}$. 
On notera cependant que les limites inductives filtrantes ne sont pas a priori représentables dans la catégorie 
$\bMod_\mQ(\bvocB)$.

\begin{prop}[\cite{ag3} 12.18]\label{higgs0-md15}
Pour tout $\bvocB_\mQ$-module de Dolbeault $\cM$, 
$\cH(\cM)$ \eqref{higgs0-md12b} est un $\co_\fX[\frac 1 p]$-fibré de Higgs soluble et associé à $\cM$.  
En particulier, $\cH$ induit un foncteur que l'on note encore
\begin{equation}\label{higgs0-md15a}
\cH\colon \bMod^\Dolb_\mQ(\bvocB)\rightarrow \bMH^\sol(\co_\fX[\frac 1 p], \xi^{-1}\Omega^1_{\fX/\cS}), 
\ \ \ \cM\mapsto \cH(\cM).
\end{equation}
\end{prop}

\begin{prop}[\cite{ag3} 12.23]\label{higgs0-md16}
On a un foncteur
\begin{equation}\label{higgs0-md16a}
\cV\colon \bMH^\sol(\co_\fX[\frac 1 p], \xi^{-1}\Omega^1_{\fX/\cS})\rightarrow \bMod^\Dolb_\mQ(\bvocB), 
\ \ \ \cN\mapsto \underset{\underset{r\in \mQ_{>0}}{\longrightarrow}}{\lim}\ \cK^{r}(\top^{r+}(\cN))).
\end{equation}
De plus, pour tout objet $\cN$ de $\bMH^\sol(\co_\fX[\frac 1 p], \xi^{-1}\Omega^1_{\fX/\cS})$, 
$\cV(\cN)$ est associé à $\cN$.
\end{prop} 

\begin{teo}[\cite{ag3} 12.26]\label{higgs0-md18}
Les foncteurs \eqref{higgs0-md15a} et \eqref{higgs0-md16a} 
\begin{equation}
\xymatrix{
{\bMod^\Dolb_\mQ(\bvocB)}\ar@<1ex>[r]^-(0.5){\cH}&{\bMH^\sol(\co_\fX[\frac 1 p], \xi^{-1}\Omega^1_{\fX/\cS})}
\ar@<1ex>[l]^-(0.5){\cV}}
\end{equation}
sont des équivalences de catégories quasi-inverses l'une de l'autre. 
\end{teo}

\begin{teo}[\cite{ag3} 12.34]\label{higgs0-md181}
Soient $\cM$ un $\bvocB_\mQ$-module de Dolbeault, $q$ un entier $\geq 0$. 
Notons $\mK^\bullet(\cH(\cM))$ le complexe de Dolbeault du $\co_\fX[\frac 1 p]$-fibré de Higgs 
$\cH(\cM)$ \eqref{higgs0-not2}. On a alors un isomorphisme canonique fonctoriel de $\co_\fX[\frac 1 p]$-modules \eqref{higgs0-md9ab}
\begin{equation}\label{higgs0-MF27a}
\rR^q\top_*(\cM)\stackrel{\sim}{\rightarrow}\rH^q(\mK^\bullet(\cH(\cM))).
\end{equation}
\end{teo}

\subsection{}\label{higgs0-md19}
Soit $g\colon X'\rightarrow X$ un morphisme étale. 
Il existe essentiellement un unique morphisme étale $\tg\colon \tX'\rightarrow \tX$ 
qui s'insère dans un diagramme cartésien \eqref{higgs0-not1a}
\begin{equation}\label{higgs0-md19a}
\xymatrix{
{\coX'}\ar[r]\ar[d]_{\cog}&{\tX'}\ar[d]^{\tg}\\
{\coX}\ar[r]&{\tX}}
\end{equation}
de sorte que $\tX'$ est une $\cA_2(\oS)$-déformation lisse de $\coX'$.
On associe à $(X',\tX')$ des objets analogues à ceux définis plus haut pour $(X,\tX)$, 
qu'on note par les mêmes symboles affectés d'un exposant $^\prime$.  
Le morphisme $g$ définit par fonctorialité un morphisme de topos annelés (\cite{ag3} 8.20)
\begin{equation}\label{higgs0-md19b}
\Phi\colon (\tE',\ocB')\rightarrow (\tE,\ocB).
\end{equation}
On démontre dans (\cite{ag3} 8.21) que $\Phi$ s'identifie au morphisme de localisation de $(\tE,\ocB)$ en $\sigma^*(X')$. 
Par ailleurs, $\Phi$ induit un morphisme de topos annelés 
\begin{equation}\label{higgs0-md19c}
\bvPhi\colon (\tE'^{\mN^\circ}_s,\bvocB')\rightarrow (\tE^{\mN^\circ}_s,\bvocB).
\end{equation} 
On note $\fgg\colon \fX'\rightarrow \fX$ 
le prolongement de $\ogg\colon \oX'\rightarrow \oX$ aux complétés $p$-adiques. 

\begin{prop}[\cite{ag3} 14.9]\label{higgs0-md20}
Les hypothèses étant celles de \eqref{higgs0-md19}, soient, de plus, $\cM$ un $\bvocB_\mQ$-module de Dolbeault, 
$\cN$ un $\co_\fX[\frac 1 p]$-fibré de Higgs soluble à coefficients dans $\xi^{-1}\Omega^1_{\fX/\cS}$.
Alors $\bvPhi^*(\cM)$ est un $\bvocB'_\mQ$-module de Dolbeault et $\fgg^*(\cN)$ est un 
$\co_{\fX'}[\frac 1 p]$-fibré de Higgs soluble à coefficients dans $\xi^{-1}\Omega^1_{\fX'/\cS}$.
Si de plus, $\cM$ et $\cN$ sont associés, $\bvPhi^*(\cM)$  et $\fgg^*(\cN)$ sont associés. 
\end{prop}

On démontre en fait que les diagrammes de foncteurs  
\begin{equation}\label{higgs0-md20a}
\xymatrix{
{\bMod_\mQ^\Dolb(\bvocB)}\ar[r]^-(0.5)\cH\ar[d]_{\bvPhi^*}&
{\bMH^\sol(\co_\fX[\frac 1 p],\xi^{-1}\Omega^1_{\fX/\cS})}\ar[d]^{\fgg^*}\ar[r]^-(0.5)\cV&{\bMod_\mQ^\Dolb(\bvocB)}\ar[d]^{\bvPhi^*}\\
{\bMod_\mQ^\Dolb(\bvocB')}\ar[r]^-(0.5){\cH'}&{\bMH^\sol(\co_{\fX'}[\frac 1 p],\xi^{-1}\Omega^1_{\fX'/\cS})}
\ar[r]^-(0.5){\cV'}&{\bMod_\mQ^\Dolb(\bvocB')}}
\end{equation}
sont commutatifs à isomorphismes canoniques près (\cite{ag3} 14.11).

\subsection{}\label{higgs0-md21}
Il existe un unique morphisme de topos 
\begin{equation}\label{higgs0-md21a}
\psi\colon \tE_s^{\mN^\circ}\rightarrow X_\et
\end{equation}
tel que pour tout $U\in \ob(\Et_{/X})$, $\psi^*(U)$ soit le système projectif constant 
$(\sigma^*_s(U_s))_{\mN}$ \eqref{higgs0-TFT8c}. 
On note $\Et_{\coh/X}$ la sous-catégorie pleine de $\Et_{/X}$ formée des schémas étales 
de présentation finie sur $X$. On a une catégorie fibrée canonique
\begin{equation}\label{higgs0-md21b}
\MOD_\mQ(\bvocB)\rightarrow \Et_{\coh/X},
\end{equation}
dont la fibre au-dessus d'un objet $U$ de $\Et_{\coh/X}$ est la catégorie $\bMod_\mQ(\bvocB|\psi^*(U))$
et le foncteur image inverse par un morphisme $U'\rightarrow U$ de $\Et_{\coh/X}$ 
est le foncteur de restriction \eqref{higgs0-md19b}
\begin{equation}\label{higgs0-md21c}
\bMod_\mQ(\bvocB|\psi^*(U))\rightarrow \bMod_\mQ(\bvocB|\psi^*(U')), \ \ \ \cM\mapsto \cM|\psi^*(U').
\end{equation}
D'après \ref{higgs0-md20}, elle induit une catégorie fibrée 
\begin{equation}\label{higgs0-md21d}
\MOD_\mQ^\Dolb(\bvocB)\rightarrow \Et_{\coh/X}
\end{equation}
dont la fibre au-dessus d'un objet $U$ de $\Et_{\coh/X}$ est la catégorie $\bMod_\mQ^\Dolb(\bvocB|\psi^*(U))$.

\begin{prop}[\cite{ag3} 15.4]\label{higgs0-md22}
Soient $\cM$ un objet de $\bMod^{\atf}_\mQ(\bvocB)$, $(U_i)_{i\in I}$ un recouvrement de $\Et_{\coh/X}$. 
Pour que $\cM$ soit de Dolbeault, il faut et il suffit que pour tout $i\in I$, le $(\bvocB|\psi^*(U_i))_\mQ$-module
$\cM|\psi^*(U_i)$ soit de Dolbeault.
\end{prop}

\begin{prop}[\cite{ag3} 15.5]\label{higgs0-md27}
Les conditions suivantes sont équivalentes~:
\begin{itemize}
\item[{\rm (i)}] La catégorie fibrée \eqref{higgs0-md21d}
\begin{equation}\label{higgs0-md23a}
\MOD_\mQ^\Dolb(\bvocB)\rightarrow \Et_{\coh/X}
\end{equation}
est un champ {\rm (\cite{giraud2} II 1.2.1)}.
\item[{\rm (ii)}] Pour tout recouvrement $(U_i\rightarrow U)_{i\in I}$ de $\Et_{\coh/X}$, 
notant $\cU$ (resp. pour tout $i\in I$, $\cU_i$) le schéma formel complété $p$-adique de $\oU$ (resp. $\oU_i$), 
pour qu'un $\co_\cU[\frac 1 p]$-fibré de Higgs $\cN$ à coefficients dans $\xi^{-1}\Omega^1_{\cU/\cS}$
soit soluble, il faut et il suffit que pour tout $i\in I$, le $\co_{\cU_i}[\frac 1 p]$-fibré de Higgs $\cN\otimes_{\co_\cU}\co_{\cU_i}$ 
à coefficients dans $\xi^{-1}\Omega^1_{\cU_i/\cS}$ soit soluble.
\end{itemize}
\end{prop}

\begin{defi}[\cite{ag3} 15.6]\label{higgs0-md23}
Soit $(\cN,\theta)$ un $\co_\fX[\frac 1 p]$-fibré de Higgs  à coefficients dans $\xi^{-1}\Omega^1_{\fX/\cS}$.
\begin{itemize}
\item[(i)] On dit que $(\cN,\theta)$ est {\em petit} s'il existe un sous-$\co_\fX$-module cohérent $\fN$ de $\cN$
qui l'engendre sur $\co_\fX[\frac 1 p]$ et un nombre rationnel $\varepsilon>\frac{1}{p-1}$ tels que 
\begin{equation}
\theta(\fN)\subset  p^\varepsilon \xi^{-1}\Omega^1_{\fX/\cS}\otimes_{\co_\fX} \fN.
\end{equation}
\item[(ii)] On dit que $(\cN,\theta)$ est {\em localement petit} s'il existe un recouvrement ouvert $(U_i)_{i\in I}$ de $X_s$
tel que pour tout $i\in I$, $(\cN|U_i,\theta|U_i)$ soit petit. 
\end{itemize}
\end{defi}

\begin{prop}[\cite{ag3} 15.8]\label{higgs0-md24}
Tout $\co_\fX[\frac 1 p]$-fibré de Higgs soluble $(\cN,\theta)$ à coefficients dans $\xi^{-1}\Omega^1_{\fX/\cS}$ est localement petit.
\end{prop}

\begin{prop}[\cite{ag3} 15.9]\label{higgs0-md25}
Supposons que $X$ soit affine et connexe et qu'il admette un $S$-morphisme étale
dans $\mG_{m,S}^d$ pour un entier $d\geq 1$. 
Alors, tout petit $\co_\fX[\frac 1 p]$-fibré de Higgs à coefficients dans $\xi^{-1}\Omega^1_{\fX/\cS}$ est soluble.
\end{prop}

\begin{cor}[\cite{ag3} 15.10]\label{higgs0-md26}
Si les conditions de \eqref{higgs0-md27} sont remplies, 
tout $\co_\fX[\frac 1 p]$-fibré de Higgs localement petit à coefficients dans $\xi^{-1}\Omega^1_{\fX/\cS}$ est soluble.
\end{cor}

\end{document}